\documentclass[10pt]{article}

\usepackage{amsfonts}
\usepackage{amsmath}
\usepackage{amssymb}
\usepackage{color}
\usepackage{epsf,psfrag}
\usepackage{latexsym}
\usepackage[pdftex]{graphicx}
\usepackage{pgfplots}
\usepackage{epstopdf}

\newtheorem{remark}{Remark}
\newtheorem{example}{Example}

\newtheorem{lemma}{Lemma}
\newtheorem{theorem}{Theorem}
\newtheorem{corollary}{Corollary}
\newtheorem{prop}{Proposition}

\title{A Critical Account of Perturbation Analysis of Markov Chains\footnote{Edited version of the paper that appeared in {\em  Markov Processes and Related Fields}, {\bf 22}, pages 227-265, 2016.}}

\author{Karim Abbas\\
 LAMOS, University of Bejaia, Algeria\\
Email: karabbas2003@yahoo.fr\\
 and\\
Joost Berkhout\\
Vrije Universiteit Amsterdam\\
Department of Econometrics and Operations Research\\
The Netherlands\\
Email: j2.berkhout@vu.nl\\
 and\\
Bernd Heidergott\\
Vrije Universiteit Amsterdam\\
Department of Econometrics and Operations Research \& Tinbergen Institute\\
The Netherlands\\
Email: b.f.heidergott@vu.nl\\}

\begin{document}

\date{ }

\maketitle

\begin{abstract}
Perturbation analysis of Markov chains provides bounds on the effect that a change in a Markov transition matrix has on the corresponding stationary distribution. This paper compares and analyzes bounds found in the literature for finite and denumerable Markov chains and introduces new bounds based on series expansions. 
We discuss a series of examples to illustrate the applicability and numerical efficiency of the various bounds.
Specifically, we address the question on how the bounds developed for finite Markov chains behave as the size of the system grows.
In addition, we provide for the first time an analysis of the relative error of these bounds. 
For the case of a scaled perturbation we show that perturbation bounds can be used to analyze stability of a stable Markov chain with respect to  
perturbation with an unstable chain. 
\\

\noindent
{\bf Keywords:} Markov chains, perturbation bounds, condition number, strong stability, series expansion, queuing \\

\noindent
AMS Primary: 60J10;  Secondary: 15A12; 15A18 
\end{abstract}

\section{Introduction}

Perturbation analysis of Markov chains (PAMC) studies the effect  a perturbation of a Markov transition matrix has on the stationary distribution of the chain.
Consider a Markov chain with discrete state space $ S $, transition probability matrix $ P $, and unique stationary distribution $ \pi_P$.
Furthermore, let  $ R $ be an alternative Markov transition matrix  on $ S $ with unique stationary distribution $ \pi_{R} $.
PAMC addresses the following question:
what is the effect of switching from $ P $ to $ R $ on the stationary distribution of the chain?
More formally,
PAMC theory studies bounds of the type
\begin{equation}\label{opl}
|| \pi_{R}^\top - \pi_P^\top|| \leq  {  \Delta} ( R , P )  ,
\end{equation}
where $ || \cdot || $ denotes a suitable vector norm (details will be provided later in the text),
$ {\Delta } ( R , P )  $ is a scalar function of $ P $ and $ R $, and $ \top $ denotes the transposed\footnote{We use the transposed here as in this paper all vectors are  by convention column vectors.}. 
The study of the effect of perturbing a Markov transition matrix on its stationary distribution dates back to Schweitzer's pioneering paper \cite{PA7}.
Best to our knowledge, the first paper putting this perturbation question into the framework of (\ref{opl}) is \cite{MeyerErst}. 
Specifically, \cite{MeyerErst} proposed  bounds of the form
\begin{equation}\label{Delta}
 {\Delta} ( R , P )  = \kappa || R -P || ,
 \end{equation}
for some appropriate matrix norm, where $ \kappa $ is the so-called {\em condition number}.
While the condition number is typically applied to bounding the effect in terms of $ R -P $, Theorem~3.2 in \cite{mitrophanov} provides
a condition number for  $ \| R^m  - P^m ||$. 
In the remainder of this article we will refer to  any instance of the  bound in (\ref{opl}) with $ { \Delta} ( R , P ) $ as in
(\ref{Delta})  as {\em condition number  bound} (CNB).

PAMC is a field of active research \cite{anisimov,kirkland,H2,neumannxu,Senata,mitrophanov,AE1,AE2,Neu1,Neu2}
and various CNBs have been proposed in the literature \cite{chomayer,AE2}. 
As we will discuss later on in more detail, an alternative type of bound called {\em strong stability bound} (SSB) can be derived via the strong stability method.
SSB bounds  the weighted supremum norm of $ \pi_{R}^\top - \pi_P^\top $ by an  expression that is non-linear as function of $ || R - P  || $. 
For early references see \cite{karta96,kartashov86} and recent references are \cite{liu12, lekadir,Rab}.
Perturbation bounds are of interest in a wide area of applications. 
For example, in mathematical physics \cite{Szehr} and climate modeling \cite{Chek}, in Bayesian statistics \cite{Andrieu,Alquier}, and in Bioinformatics \cite{AY2, Pal}.
Perturbation bounds have also been applied in robustness analysis of social networks and of Google's PageRank algorithm \cite{extra}.


A fruitful model for PAMC is that of a scaled perturbation.
More specifically, let $ R , P $ be two Markov kernels defined
on the same state space. Then the convex combination of both kernels
\begin{equation}\label{eq:convex}
P ( \theta ) = ( 1 - \theta ) P + \theta R , \quad \theta \in [ 0 ,1 ] ,
\end{equation}
 is a well-defined Markov kernel. 
Note that $ P ( 0 ) = P $ and $ P ( 1 ) = R $.
In perturbation analysis of $P(\theta)$ we are interested in the effect of changing $ \theta  $ from $ 0 $ to some 
value $0 < \theta  \leq 1$.
By linearity of norms,
\begin{equation}\label{eq:deltanorm}
 || P ( \theta)  - P   || = \theta || R - P || ,
\end{equation}
for $\theta \in [0,1]$.
This allows to scale the size of the perturbation via control parameter $ \theta$.
Letting  
\begin{equation}\label{relerr}
\eta ( R , P ) = \frac{  \Delta ( R , P ) -  || \pi_R^\top - \pi_P^\top ||}{  || \pi_R^\top - \pi_P^\top ||}
\end{equation}
denote the {\em relative error of perturbation bound $  {\Delta} ( R , P ) $},
scaled perturbations, i.e., $ R = P ( \theta ) $, allow for analyzing the
behavior  of the relative error $ \eta ( \theta ) = \eta ( P ( \theta ) , P ) $ as $ \theta $ tends to zero. 

The analysis of scaled perturbation is of particular  interest if $ P ( \theta ) $,  for $ \theta \in [ 0 , 1] $, has a clear interpretation. We will illustrate this by a queueing model with denumerable state-space and breakdowns, where $ \theta $ models the probability of a breakdown.
An interesting observation is that in the parametrized model we establish conditions for stability of a mixture of a stable (no breakdowns) and an unstable (only breakdowns) Markov chain modeling a pure birth process.
More specifically, we apply PAMC techniques to provide a lower bound for the domain of stability of $ P ( \theta ) $.
The contributions of the paper are the following: 
\begin{itemize}

\item We provide a unified approach to PAMC for finite and denumerable Markov chains. Our analysis covers CNBs and SSB.  

\item We introduce new bounds that do have the desirable property  that the relative error of the bound tends to zero as the 
size of the perturbation tends to zero. 
These new bounds are derived by a series expansion approach.   

\item We will provide sufficient conditions under which the convergence of the series expansion already constitutes 
existence of a stationary distribution. By introducing the new concept of the bias term, we are able to treat the case of Markov multi-chains (i.e., chains with several ergodic classes) and uni-chains in a unified framework.

\item We will show that techniques derived in PAMC can be applied to stability analysis. 
A worked out example from queuing theory will illustrate the fruitfulness of PAMC methods for this type of problem.

\end{itemize}


The paper is organized as follows.
In Section~\ref{sec:3}  the perturbation bounds are presented and the main theoretical results a established.
For a simple example, Section~\ref{sec:leadex} presents explicit solutions for the various bounds.
Section~\ref{sec:XXX} is devoted to perturbation bounds for the  M/G/1 queue with breakdowns.
Other than the small numerical examples reported in the literature, 
the queuing system will be analyzed for the case of a large but finite state-space and for the infinite dimensional case.

\section{Perturbation Analysis}\label{sec:3}

Throughout this paper we will consider aperiodic Markov chains defined on an at most denumerable state space $ S  = \{ 0 , 1, \ldots \} \subset \mathbb{N}$.
Unless stated otherwise, we assume that the Markov chains are aperiodic and have one closed communicating class of states with possible transient states.

\subsection{Preliminaries and Basic Definitions}

If $ P = (P_{ij})_{i,j\in S}$ is a Markov transition matrix of some Markov chain $ \{X_k \} $, then
$ P_{ i j } = \mathbb{E} [ 1_{ X_{k+1} } ( j ) | X_k = i ] $ for $ i , j \in S $ and $ k\in\mathbb{N}$,  
where $ 1_j (i) $ is one if $ j=i $ and zero otherwise, $i, j \in S$. Sometimes $P(i,j):=P_{ij}$ is used instead for notational clarity. Further, let $f \in \mathbb{R}^S$ be a reward vector where $f_i$ is the reward for being in state $i\in S$. 
With these definitions, one obtains
\begin{equation}\label{eq:tr}
 \mu^\top P f = \sum_{ i , j  \in S } \mu_i \! P _{ i j } f_j =  \sum_{i \in S } \mathbb{E} [  f_{ X_1}  | X_0 = i ] \mu_i
\end{equation}
as the expected reward after one transition provided the Markov chain is started in state $ i $ with
probability $ \mu_i$, for $i \in S$, in vector-matrix notation\footnote{As exemplified in  (\ref{eq:tr}), 
distributions on $ S $ are represented  as row vectors in vector-matrix notation.
Since by convention a tuple $ \mu= ( \mu_i : i \in S ) $ representing a distribution on $ S $ 
when written as vector $ \mu \in [0,1]^S$ becomes a column vector, 
we explicitly denote $ \mu $ in transposed form to make it a row vector, i.e., we write $ \mu^\top$,
see  (\ref{eq:tr}).
When it causes no confusion we will refer to either $ \mu $ or $ \mu^\top $ as distributions.
For example, in (\ref{eq:tr}) we may refer to $ \mu $ as well as $ \mu^\top $ as initial distribution.}.
For more details we refer to \cite{heder4-11,sampledchain}.

In the following denote the ergodic projector of $ P $ by $\Pi_P$, i.e., the matrix with rows identical to $ \pi_P^\top $, and we let
$ D_P  $ denote the deviation matrix of $ P $, which is given by
\begin{equation}\label{eq:ddd}
D_P=  \sum_{k=0}^\infty ( P^k - \Pi_P) = ( I - P  + \Pi_P ) ^{-1} -  \Pi_P  ,
\end{equation}
provided that it exists.
The matrix $ ( I - P  + \Pi_P ) ^{-1} $ is called the {\em fundamental matrix} (potential) of $ P $.
Letting $ A^{\#} $ denote the group inverse of the matrix $ A = I - P $, see \cite{a,bbb}, it holds that
$ D_P = A^{\#}$ if the deviation matrix exists. 
Conditions for existence of the deviation matrix and its related properties have been extensively studied in the
literature, see \cite{heder4-11,Syski}.
For finite Markov chains, the deviation matrix is an instance of the generalized inverse of $ I - P $;
see \cite{bbb} for an early reference. As Hunter demonstrates in \cite{2} for finite Markov chains, the generalized inverse plays a major role
in perturbation analysis.

For $ x \in \mathbb{R}^S $, we denote by $ || x ||_\infty $ the maximum absolute value (a.k.a.~infinity norm or $\infty$-norm) and by $ || x ||_1 $ the sum of absolute values (a.k.a.~$L_1$ norm or 1-norm). 
Furthermore, for  $ v $, such that $ v (i) \geq 1 $ for all $ i \in S $ and $ v ( 0 ) =1 $, we denote by
\begin{equation}\label{eq:normf}
 \| x  \|_{\upsilon} = \sup_{i \in S} \frac{|x_i |}{\upsilon(i)}
\end{equation}
the weighted supremum norm of $ x\in \mathbb{R}^S$, also called {\em  $v$-norm}.
In the following we let
\[
 v_\alpha ( i ) = \alpha^i , \quad  i \in S ,
\]
with $ \alpha \in [ 1 , \infty ) $ some unspecified constant.
In the following we will omit  the subscript $ \alpha  $ whenever the results stated hold for general $ \alpha \geq 1$.
Norms are extended to matrices by using the standard induced norms and 
vector norms are obtained as the vector norm induced by the corresponding matrix norm\footnote{Note that this implies for $ x \in S^\mathbb{R}$:  $ || x^\top ||_\infty = || x ||_1 $ and $|| x^\top ||_ 1 = \| x \|_\infty $.}
Note that as probability measures are row vectors, the $ v$-norm of a measure $ \mu^\top $ on $ S $ is given by 
\begin{equation}\label{eq:normmu}
\|  \mu^\top \|_{\upsilon} = \sum_{k \in S} v ( k ) | \mu_{k}| .
\end{equation}
Specifically, applying the $ v $-norm to (\ref{eq:tr}) one readily obtains
\[
|  \mu^\top  P f  | \leq  || \mu^\top ||_v \, || P ||_v  \, || f ||_v ,
\]
which shows that  (\ref{eq:normmu}) and  (\ref{eq:normf}) arise naturally in PAMC.
Note that by (\ref{eq:normmu})  for a (possibly signed) measure  $\mu  $ on $ \mathbb{R}^{S}$
the $v$-norm of $ \mu^\top $, for $ v \equiv 1$, coincides with total variational norm.
Note further  that $ v \geq 1 $ implies for $ x \in \mathbb{R}^S $ that  $ \|  x\|_v \leq \| x \|_\infty $.
In addition, for a measure $ \mu$ we have that $ \| \mu^\top \|_{1 } \leq \| \mu^\top \|_{\infty } \leq \| \mu^\top  \|_v $.

In the following we will omit the norm type indicator and use the generic $ || \cdot || $ sign whenever the result holds for any of the above norms.
If a result is limited to a particular norm, this will be clearly indicated.

To illustrate the efficiency of the bounds we will use throughout the paper three different types of finite Markov chains introduced in the following. An example of a Markov chain on a denumerable state space will be discussed in detail in 
Section~\ref{sec:A5}.

\begin{example}\label{ex:mc} 
{\bf Two-State Chain:} Let $S = \{ 0, 1\}$ and 
\[
P^s = \begin{pmatrix}
      1 - p & p \\
      q & 1 - q \\
    \end{pmatrix},
\]
with $p, q \in (0, 1)$. It is easily checked  that 
\[
\pi_{P^s} = \frac{1}{ p + q} (q , p)^\top
\]
is the stationary distribution of $ P^s $. 
The deviation matrix is given by: 
\[
D_{P^s} = \frac{ 1 }{ ( p + q )^{2} } \begin{pmatrix}
                                       p & - p \\
                                       - q & q \\
                                     \end{pmatrix}.
\]

\noindent
{\bf Ring Network:} The next example that we will discuss is that of a ring, introduced in the following.
Let $S = \{0, \ldots , n-1 \}$ and for any $n \geq 2$,
\[
P^\circ (n) =
\left(
              \begin{array}{cccccc}
                1- 2 b   & b  & 0 & 0 & \ldots & b \\
                b  & 1 - 2 \, b  & b & 0 & \ldots & 0 \\
                0 & b & 1 - 2 \, b  & b & \ldots & 0 \\
                \vdots & \vdots & \ddots & \ddots & \ddots & \vdots \\
                0 & 0 & \ldots & b & 1 - 2 \, b & b \\
                b  & 0 & \ldots & 0 & b & 1 -  2 b  \\
              \end{array}
            \right),
\]
with $b \in (0, 1/2]$. We get the stationary distribution:
\[
\pi^\circ_i (n) = \frac{1}{n}, \;\; \mbox{ for } i \in S.
\]
For the deviation matrix, we obtain:
\[
D^\circ(n) :  = D_{ P^\circ (n) } = \left(
                \begin{array}{ccccc}
                  d_{0} & d_{1} & d_{2} & \ldots & d_{n-1} \\
                  d_{n-1} & d_{0} & d_{1} & \ldots & d_{n - 2} \\
                  d_{n - 2} & d_{n-1} & d_{0} & \ldots & d_{n - 3} \\
                  \vdots & \vdots & \vdots & \ddots & \vdots \\
                  d_{1} & d_{2} & d_{3} & \ldots & d_{0} \\
                \end{array}
              \right),
\]
where
\begin{equation}\label{di:eq}
    d_{i} = \frac{ (n - 1) (n + 1) }{12 \, b \, n} - \frac{ (n - i) i  }{2 \, b \, n} \, \, \, \hbox{for} \, \, \, i \in S.
\end{equation}
Furthermore, $\sum_{i = 0}^{n-1} d_{i} = 0$. Equivalently, $D^\circ (n)$ can be expressed as
\[
D^\circ (n) = \left( \widetilde D_{ij}(n) \right)_{i,j \in S},
\]
where
\begin{equation}\label{dij:eq}
\widetilde D_{ij} (n) = d_{(j - i)(\mbox{\emph{mod}} \; n) + 1} = \frac{ (n - 1) (n + 1) }{12 \, b \, n} -\frac{ \{ n - (j - i) (\mbox{\emph{mod}} \; n)\} \{ (j - i) (\mbox{\emph{mod}} \; n)\}}{2 \, b \, n}.
\end{equation}

\noindent
{\bf Star Network:}
The third example considered is the Star Network with state space $S=\{0,\dots,n-1\}$. For $n\geq 2$ let
\[
P^\star (n)=
\left(
              \begin{array}{cccccc}
               1-  \beta  &  \frac{ \beta }{ n - 1 } &  \frac{ \beta }{ n - 1 } &  \frac{ \beta }{ n - 1 } & \ldots &  \frac{ \beta }{ n - 1 } \\
               1 - \gamma &  \gamma & 0 & 0 & \ldots & 0 \\
                1- \gamma  &  0 &  \gamma & 0 & \ldots & 0 \\
                \vdots & \vdots & \ddots & \ddots & \ddots & \vdots \\
                1- \gamma   & 0 & \ldots & 0 &   \gamma & 0 \\
                1  -\gamma & 0 & \ldots & 0 & 0 &  \gamma \\
              \end{array}
            \right) ,
\]
for $\beta \in (0, 1]$ and $\gamma \in [0, 1)$.
Following \cite{golub2010}, the stationary distribution is given by
\[
\pi_i^\star(n) =
\left\{
               \begin{array}{ll}
                 \frac{ 1 - \gamma }{ 1 - \gamma + \beta } & \hbox{ for } \, i = 0, \\
\\
                 \frac{ \beta }{ ( n - 1 ) ( 1 - \gamma + \beta )} & \hbox{ for } i > 0.
               \end{array}
             \right.
\]
For the deviation matrix, we obtain:
\[
D^\star(n) = \left(
                   \begin{array}{ccc}
                     \frac{ \beta }{ ( 1 - \gamma + \beta )^{ 2 } } & \vline & - \frac{ \beta }{ ( 1 - \gamma + \beta )^{ 2 } ( n - 1 ) } \bar{1}^{T} \\
               ------   & \vline &  ------------- \\
                     - \frac{ ( 1 - \gamma ) }{ ( 1 - \gamma + \beta )^{ 2 } } \bar{1} & \vline & \frac{ 1 }{ ( 1 - \gamma ) } I - \frac{ \beta \{ ( 1 - \gamma ) + ( 1 - \gamma + \beta ) \} }{ ( 1 - \gamma )  ( 1 - \gamma + \beta )^{ 2 } ( n - 1 ) } \bar{1}\bar{1}^\top\\
                   \end{array}
                 \right),
\]

where $\bar{1} = [1, \ldots, 1]^{T}$ of size $n-1$ and $I$ denotes the $(n-1)\times(n-1)$ identity matrix.
\end{example}

In our analysis we will frequently work with the {\em taboo kernel} of a Markov transition matrix $ P $.
In \cite{kartashov86}  a very elegant and flexible way for obtaining a taboo kernel is described.
For this let $ h $ be a non-negative vector and $ \sigma $  a probability measure on $ S $, 
such that $\pi^\top_P h > 0$ and $ P -  h \sigma^\top  $ is matrix with non-negative values, where $ h \sigma^\top  $ denotes the matrix product of vector $ h $ and $ \sigma $, i.e., $ h \sigma^\top $ is a square matrix.
Then, the taboo kernel of $ P $ with respect to $ h $ and $ \sigma $ is defined as
\begin{equation}\label{eq:taboo}
T := P  - h   \sigma^\top .
\end{equation}
For example, let
\[
h   = ( P (0 ,0 ), P (1 , 0), P ( 2 , 0) , ...)^\top
\]
denote the first column of $ P$, and let $ \sigma = ( 1 , 0 , 0 , \ldots )^\top $,
then
\[
 T = P - h \sigma^\top =
\left \{ \begin{array}{cc}  P (
i , j ) & j > 0 \\
0 & \text{ otherwise }
\end{array} \right. .
\]
In words, $  T $ is a degenerate transition kernel that avoids entering
state zero which is obtained by setting the first column of $ P $ to zero.
Alternatively, letting $ h = ( 1 , 0 , 0 , \ldots )^\top $ and
\[
\sigma  = ( P (0 ,0 ), P (0 , 1), P ( 0 , 2) , ...)^\top ,
\]
then $ T =  P - h \sigma^\top $ is a degenerate transition kernel that never leaves state zero,
which is obtained by setting the first row of $ P$ to zero.
The taboo kernel is also known as the {\em residual matrix} in the literature, see \cite{MA}.
  
In the following we write $ _i \! P $ for the  degenerate transition kernel that avoids entering
state $i$  which is obtained by setting the $i$th column of $ P $ to zero, i.e., letting
$  \sigma = ( 0 , \ldots 0 , 1 , 0 , \ldots  )^\top $, where the entry 1 is at the $ i$th position, and $ h $ is  the $ i$th column of $ P$.
The taboo kernel $ _i \! P $ provides a convenient sufficient condition for positive recurrence of $ P $ on a denumerable state space.  The precise statement is provided in the
following proposition.

\begin{prop}\label{cor:pos}
Let $ P $ be irreducible. If for at least one $ i \in S $ it holds that $ || _i \! P || < 1 $, then
$ P $ is positive recurrent.
\end{prop}
{\bf Proof:}
First note that the $(j,k)$-th element of $ \sum_{n=0}^\infty ( _i P )^n $ gives the expected number of visits to state $k$ before jumping to state $i$ when starting in state $j$. The mean recurrence time at state $ i $ is thus given by summing the $i$-th row of $ \sum_{n=0}^\infty ( _i P )^n$, which is finite due to the norm condition.
Therefore, state $ i $ is positive recurrent. From irreducibility of $ P $ it follows that all states are positive recurrent.
\hfill $ \Box $

We call $  T $ \emph{proper} if $ || T || < 1$.
Provided that $ T  $ defined in (\ref{eq:taboo}) is proper in case of the $v$-norm, the $ v$-norm of $ \pi_P^\top $ can be bounded by
\begin{equation}\label{eq:taboo2}
\| \pi_P^\top \|_{v} \leq \frac{  \pi_P^\top  h \| \sigma^\top \|_v  }{ 1 - ||  T ||_v },
\end{equation}
see \cite{kartashov86}.
Moreover, if $ T $ is proper, the deviation matrix can alternatively be written as
\begin{equation}\label{eq:D}
D_P = ( I - \Pi_P )   \sum_{n=0}^\infty T^n  ( I -  \Pi_P )  =  ( I - \Pi_P ) ( I - T)^{-1}  ( I -  \Pi_P ) ,
\end{equation}
see \cite{H1,karta96}, where used that 
\begin{equation}\label{eq:hordijk}
( I - T)^{-1}  =  \sum_{n=0}^\infty T^n   ,
\end{equation}
provided that $ || T || < 1 $.

The idea behind considering $ T $ rather than $ P $, is that $  T $ might be constructed in such a way that the norm of $ T  $ is strictly less than one.
The following example illustrates the effect on $ || T ||$ from either removing the first column or first row. Note that removing the second column or second row may  lead to other values of $|| T ||$.

\begin{example}\label{ex:T}
For the two-state chain, i.e., $P=P^s$, we find after removing the first column
\[
  || T ||_{v} = \max\{ \alpha p, 1 - q \}.
\]
Removing the first row leads to 
\[
  || T ||_{v} =  \frac{(1-\alpha)q}{\alpha} + 1  .
\]
For the Ring and the Star networks we present the resulting norms for $||T||_{v}$ (including the $\infty$-norm  by letting $ \alpha $ tend to 1) and $||T||_{1}$, respectively, in Table~\ref{tab:1} and Table~\ref{tab:2}. 

\begin{table}[ht!]
\begin{center}
\begin{tabular}{|c|c|c|}
  \hline
  Removing: & Ring (i.e., $P=P^\circ (n)$) & Star (i.e., $P=P^\star (n)$)  \\	\hline
 $1$st row of $P$  & $ \frac{ b }{ \alpha } + 1 - 2 \, b + \alpha \, b $ & $ \frac{ b }{ \alpha } + 1 - 2 \, b + \alpha \, b $  \\
 $1$st column of $P$  & $ \max \{\alpha \, b + \alpha^{n - 1} \, b, \frac{ b }{ \alpha } + 1 - 2 \, b + b \, \alpha \}$  &   $\max\left\{ \gamma,  \frac{ \alpha \, \beta }{ n - 1}\frac{1 - \alpha^{n - 1}}{ 1 - \alpha } \right\}$   \\
  \hline
\end{tabular}
\caption{The $v$-norm for different choices for $T$ (including the $\infty$-norm). \label{tab:1}}
\end{center}
\end{table}
\begin{table}[ht!]
\begin{center}
\begin{tabular}{|c|c|c|}
  \hline
  Removing: & Ring (i.e., $P=P^\circ (n)$) & Star (i.e., $P=P^\star (n)$)  \\	\hline
   $1$st row of $P$ & $ 1 $ & $ \max\{\gamma,(n-1)(1-\gamma)\} $  \\
 $1$st column of $P$  & $ 1 $  & $\gamma + \frac{\beta}{n-1}$  \\
  \hline
\end{tabular}
\end{center}
\caption{The 1-norm for different choices for $T$.\label{tab:2}}
\end{table}
\end{example}

In the following we discus a general way of choosing $ T $.
Let $ P_{ \bullet j } $ denote the $ j$-th column of $ P $.
For a column vector $ x$ we let $ \| x \|_{\inf} = \inf_i | x_i | $.
We denote the $ j$th unit vector by $ e_j $, i.e., $ e_j $ has all elements zero except for the $j$-th element which is equal to 1.

\begin{lemma}\label{le:taboo}
Let $ P $ be a Markov transition matrix on $ S $.
Let $ j^\ast $ be the column index with maximal  value  $ \| P_{ \bullet j } \|_{\inf}$.
If $  ||   P_{  \bullet  j^\ast }  ||_{\inf } > 0  $, let $ h = P_{ \bullet j^\ast} $ and $ \sigma =e_{j^\ast}$, then  for $ T $ defined as in
(\ref{eq:taboo}) it holds that $ || T ||_v < 1  $, where $ v \equiv 1 $.
\end{lemma}
{\bf Proof:}
Without loss of generality assume that after appropriate relabeling of the states $ j^\ast = 0 $.
Let  $||  P_{  \bullet 0 }  ||_{\inf } = q> 0 $.
Removing the first column from $ P $ thus decreases the row sum of each row of $ P $ by at least $ q$, which implies
the desired result.
\hfill $ \Box $ 

\subsection{Condition Number Perturbation Bounds for Finite Chains}\label{sec:norms2}

Several condition numbers have been proposed in the literature for finite Markov chains with state space $ S = \{0,1,\dots,n-1\}$, see  \cite{chomayer} for an overview.
We keep the numbering as in  \cite{chomayer}, where seven different condition numbers were discussed.
Moreover it is shown in  \cite{chomayer} that condition numbers $ \kappa_3 $ and $ \kappa_ 6$, to be defined presently, outperform the other condition numbers,
while the choice between   $ \kappa_3 $ and $ \kappa_ 6$ depends on the choice of norms.
Condition number $\kappa_{3}$ is given by \cite{haviv,Kirk}:
  \begin{equation}\label{kappa:3}
\kappa_{3} (P)=\frac{ \max_{j}( D_P(j,j) - \min_{i} D_P(i,j) ) }{ 2 }
\end{equation}
 and leads to the following bound:
\[
 || {\pi}_R^\top - \pi_P^\top ||_{1}  \leq  \kappa_{3}(P)  || R  - P ||_{\infty} .
\]
Alternatively,  condition number $\kappa_{6}$ in \cite{Senata} is given by:
\begin{equation}\label{kappa:7}
    \kappa_{6}(P) =  \frac{1}{2}  \max _{ i , j } \sum_{k=0}^{n-1} | D_P(i,k) - D_P(j,k) | ,
\end{equation}
and the resulting bound is as follows:
\[
 || {\pi}_R^\top - \pi_P^\top ||_{\infty}  \leq  \kappa_{6}(P)  || R  - P ||_{\infty} .
\]

\begin{example}
The condition numbers for the Markov chains introduced in Example~\ref{ex:mc} are as follows:
\[
  \kappa_{3} ( P^s )  =  \frac{1}{2(p+q)}
\quad \mbox{and} \quad
 \kappa_{6} ( P^s ) = \frac{1}{p+q} ,
\]
\[
  \kappa_{3} ( P^\circ (n)  )  =  \frac{\lfloor \frac{ n }{ 2 }\rfloor (n - \lfloor \frac{ n }{ 2 }\rfloor ) }{ 4 \, b \, n} ,
 \]
 \[
  \kappa_{6} ( P^\circ (n)  )  =  \frac{1}{2}\sum\limits_{k=0}^{n-1} \left |D_{P^\circ (n)}
  \left ( \left \lfloor \frac{ n }{ 2 } \right \rfloor + 1, k-1 \right )-D_{P^\circ (n)}(1, k-1) \right | ,
\]
and
\[
  \kappa_{3} ( P^\star (n))   =  \frac{ 1 }{ 2 ( 1 - \gamma ) } 
\quad \mbox{and} \quad
 \kappa_{6} ( P^\star (n))   = \frac{1}{1-\gamma}.
\]

It is worth noting that $  \kappa_{3} ( P^\circ (n)  ) $ grows linearly in $ n$. As the condition number applies to the 1-norm of
 $ \pi_R^\top - \pi_P^\top $, which is bounded by 1,
the bound becomes thus trivial for large $ n $. 
For the Star Network, $ \kappa_3 $ and $ \kappa_6 $ do not depend on $ n$ but become trivial for $ \gamma $  close to $ 1 $.

The fact that $ \kappa_3 $ and $ \kappa_6 $ behave so different for the Ring Network and the Star Network stems from the fact that both condition numbers are defined via the deviation matrix. 
The elements of the deviation matrix are related to mean recurrence times of the corresponding Markov chain, see  \cite{bbb,2}.
Specifically, in the Ring Network the length of a path from, say,  node 0 to node $ \lfloor n/2 \rfloor $ grows with $ n $, 
whereas in the Star Network any node can reached from any other node in 2 steps.
\end{example}

It is known that $ \kappa_3 ( P ) < \kappa_6 ( P ) $ (in fact it holds that $2\kappa_3(P) \leq \kappa_6(P)$), see \cite{Kirk}). Note that this inequality implies for the Ring Network that $ \kappa_6 (P^\circ (n)) $ tends  to infinity as well.
In \cite{Kirk} it is shown that $ \kappa_3 (P) \geq (n-1)/(2n)$, with $ n$ being the size of transition matrix, and a Markov chain is provided for which
equality is reached.
As we will discuss in the subsequent section, $ \kappa_6 ( P ) $ may be preferable to $ \kappa_3 ( P ) $ in case
bounds on perturbations of expected rewards are considered.
\subsection{The Choice of Norms in Perturbation Analysis}\label{sec:norms}

In bounding perturbations it is important to understand how a perturbation of the Markov chain affects the steady-state reward.
Put differently, using the notation as already introduced in the introduction, relating a perturbation bound
for $ || \pi_R^\top - \pi_P^\top ||   $ to that of $ | \pi_R^\top  f - \pi_ P^\top  f  | $ is of importance in applications. 
The following lemma formalizes how the steady-state reward can be bounded via perturbation bounds for $ || \pi_R^\top - \pi_P^\top ||   $ in case of different norms.

\begin{lemma}\label{le:norms}
For arbitrary  measures  $ \widetilde{\mu} $ and $  \mu $ on $ \mathbb{R}^S$ and cost function $ f \in  \mathbb{R}^S $  such that $ | \widetilde{\mu}^\top  f - \mu^\top f | < \infty $ it holds 
\[
| \widetilde{\mu}^\top  f - \mu^\top f |
\leq  \left \{
\begin{array}{c}
||  \widetilde{\mu}^\top  - \mu^\top  ||_\infty \, || f ||_\infty  \\
 \\
||  \widetilde{\mu}^\top  - \mu^\top  ||_1 \, || f ||_1  \\
\\
||  \widetilde{\mu}^\top  - \mu^\top  ||_v \, || f ||_v \\
\end{array} \right  . .
\]
\end{lemma}
{\bf Proof:}
By simple algebra,
\[
| \widetilde{\mu}^\top f - \mu^\top f | \leq \sum_i |  \widetilde{\mu}_i - \mu_i|  \, | f_i |
\leq \sup_j | f_j |  \sum_i |  \widetilde{\mu}_i - \mu_i|   = ||  \widetilde{\mu }^\top  - \mu^\top  ||_\infty \, || f ||_\infty  .
\]
For the last  inequality, which coincides with the second inequality in case of $\alpha = 1$, note that
\begin{align}
| \widetilde{\mu}^\top f - \mu^\top f |  &  \leq \sum_i |  \widetilde{\mu}_i - \mu_i|  \, | f_i |  \\
& =  \sum_i |  \widetilde{\mu}_i - \mu_i|  v_i  \,  \frac{ | f_i|}{v_i} \\
& \leq  \left ( \sup_j \frac{| f_j |  }{v_j } \right )  \sum_i |  \widetilde{\mu}_i - \mu_i|  v_i
\\
& =  ||  \widetilde{\mu}^\top  - \mu^\top  ||_v \, || f ||_v ,
\end{align}
which concludes the proof. 
\hfill $ \Box $ 

In this chapter we study the case that $ \mu $ in Lemma~\ref{le:norms} is a stationary distribution.
Lemma~\ref{le:norms} illustrates that there is a trade-off in the choice of norms.
Indeed, since $ \|  \pi_R^\top  - \pi_P^\top  \|_1 \leq  \|  \pi_R^\top  - \pi_P^\top    \|_\infty $ it seems attractive to
ask for perturbation bounds on $ \|  \pi_R^\top - \pi_P^\top  \|_1  $.
The downside is that this choice affects the norm of the reward vector, in particular, it holds that $\|f\|_\infty \leq \|f\|_1$.
As an illustration,  consider the following example of a finite Markov chain.
Let $ P $ be the transition matrix of a M/M/1/N queue, where $ N $ is the size of the buffer of the queue including the service place, and suppose that we are interested in the effect that replacing $ P $ by $ R $ has on the stationary queue length. 
More specifically, let $ f_l( s ) = s $, for $ s \in S = \{ 0 , 1 , \ldots , N \} $, and note that 
\[
\| f_l \|_1 =  \frac{N ( N+1)}{2} > N = \| f_l \|_\infty .
\]
In the light of Lemma~\ref{le:norms},  in bounding 
$ | \pi_P^\top f_l - \pi_R^\top f_l |$
the smaller bound on the norm distance of $ \pi_R^\top  - \pi_P^\top   $ by applying the 1-norm might be outweighed by the increase in norm for the reward.
If, on the other side, one is only interested in an overflow probability, i.e., $ f_p ( s ) = 0 $ for $s < N$ and $ f_p (N) = 1$, then
$ \| f_p \|_1 = \| f_p \|_\infty=1 $ and the 1-norm bound for $ \pi_R^\top  - \pi_P^\top   $ is appropriate.
Another example where this norm trade-off is relevant is in the analysis of the  `wisdom of crowds' phenomenon in  social networks, \cite{golub2010}.
Here, $ f $ represents a belief vector with bounded support, i.e., $ f ( s) \in [ a , b ] $ for $ a < b \in \mathbb{R} $, and
$ \pi_P^\top f $ is the consensus reached in the social network modelled by $ P $.
From the above discussion it is clear that the choice of the norm for evaluating $ \pi_R^\top - \pi_P^\top $ depends on the application.
 
In the light of the above discussion it is worth noting that the $v$-norm can be adjusted to the problem under consideration.
To see this, recall that we have assumed that $v $ is of the form $ v ( i ) = \alpha^i $, $ i \in S $, with $ \alpha $ some unspecified constant.
Let us express this dependency of $ v $ on $ \alpha $ here by writing $ v_\alpha $.
Hence, the best bound for $ | \widetilde{\mu}^\top f - \mu^\top f | $  by means of the $ v $-norm is given by the solution of the following minimization problem
\begin{equation}\label{eq:starrr}
| \widetilde{\mu}^\top f - \mu^\top f | \leq \min_\alpha ||  \widetilde{\mu}^\top  - \mu^\top  ||_{v_\alpha}  \, || f ||_{v_\alpha}  .
\end{equation}
The upside of this minimization is that it trades off the effect the norm has on the reward and the measure distance.
The downside is of course that the minimization itself  can be rather demanding as $ ||  \widetilde{\mu}^\top  - \mu^\top  ||_{v_\alpha} $ or a bound thereof typically
has a complex form.
For denumerable Markov chains, $ v $ can be constructed via a Lyapunov-type of drift condition; see \cite{liu12} for details.

\subsection{Perturbation Bounds}

In perturbation analysis, $ D_P $ occurs in conjunction with a perturbation matrix $ \Delta = R -P $ which has row sums equal to zero.
From $ \Delta ( I - \Pi_P ) = \Delta $ and \eqref{eq:D} it follows that
\begin{equation}
\label{eq:DnewState}
\Delta ( I - T)^{-1}  ( I - \Pi_P )  = \Delta D_P
\end{equation}
and instead of $ D_P $ for perturbation bounds it suffices to consider
\begin{equation}\label{eq:Dnew}
( I - T)^{-1}  ( I - \Pi_P ) .
\end{equation}
Note that due to the fact that $ \Delta ( I - T)^{-1} $ fails to have row sums equal to zero, the term $  I - \Pi_P $ on the LHS in \eqref{eq:DnewState} cannot be disregarded. In other words,
$  \Delta ( I - T)^{-1}  \not = \Delta   ( I - T)^{-1}   ( I - \Pi_P  )   $, except for special cases.
By simple algebra, it holds for Markov transition matrices $ R  $ and $ P $ that
\begin{eqnarray}
\pi_{ R }^\top & = &  \pi_{P }^\top  + \pi_{R}^\top ( R - P) D_{ P }  \label{updateD} \\
&  =  &
\pi_{P }^\top + \pi_{R}^\top ( R - P )   ( I - T)^{-1}   ( I - \Pi_P  )    \label{update} .
\end{eqnarray}

\begin{remark}\label{re:1}
The above formula is called {\em update formula} and allows for deriving a first perturbation bound.
Using the fact that $ || \pi_{R}^\top ||_\infty = 1 $, (\ref{update}) yields
\[
|| \pi_{ R}^\top  -  \pi_{P }^\top  ||_\infty \leq   ||  R- P ||_\infty \, ||  ( I - T)^{-1}  ( I - \Pi_P )  ||_\infty ,
\]
which provides a first perturbation bound.
Put differently $   ||  ( I - T)^{-1}  ( I - \Pi_P )  ||_\infty $ yields a condition number for $|| \pi_{ R}^\top  -  \pi_{P }^\top  ||_\infty$.
\end{remark}

Repeated insertion of the expression for $ \pi_{ R} $ in  (\ref{updateD}) on RHS of  (\ref{updateD}), yields
\begin{equation}\label{eq:21}
\pi_{ R }^\top = 
\pi_{P}^\top  \sum_{ k=0}^N
( ( R - P )  D_P )^k     +  \pi_{R} ^\top( ( R - P ) D_P )^{N+1}  .
\end{equation}
We call 
\[
B ( R , P ) = \lim_{N \rightarrow \infty }  \pi_{R}^\top ( ( R - P ) D_P )^{N}
\]
the {\em bias term}, provided that the limit exists.
Letting $  N $ tend to infinity in (\ref{eq:21}) 
we arrive at 
\begin{eqnarray}
\pi_{ R }^\top  
& = & 
\pi_{P}^\top  \sum_{ k=0}^\infty
( ( R - P )  D_P )^k  \: + 
B ( R , P )  \label{eq:basicseries}
\\
& = & \pi_{P}^\top  ( I - ( R - P ) D_P )^{-1}  \nonumber
+ 
B ( R , P ) ,
\end{eqnarray}
provided the series exists and the bias term is finite.
As we will explain in the following, the bias term is typically zero in case that $ R $ and $ P $ are uni-chain.
The series in (\ref{eq:basicseries}) already appears without the bias term in \cite{PA7}.
It has been rediscovered in \cite{Cao} and extended to Markov chains on a general state-space in \cite{heidergott}, both
references study problem classes where the bias term is zero.

In deriving the series expansion in (\ref{eq:basicseries}) we required that the stationary distribution $ \pi_{R } $ exists.
As the next theorem shows, convergence of the series already implies existence of $ \pi_{R } $.
Moreover, sufficient conditions are provided for the bias term to be equal to the zeros vector.

\begin{theorem}\label{th:stability}
Let $ P $ be irreducible, aperiodic and positive recurrent.
Suppose that the series in (\ref{eq:basicseries}) converges to some finite limit $ \mu^\top $, i.e., let
\[
\mu^\top = \pi_{P}^\top  ( I - ( R - P ) D_P )^{-1} .
\]
\begin{itemize}

\item[(i)]

If $ \mu_i \geq 0 $, for $ i \in S $, then $ \mu $ is a stationary distribution of $ R $.

\item[(ii)] 

If $ R $ is irreducible and aperiodic and there exists $ i \in S $ such that $ || _i R  || < 1 $, then $ \mu $ is the unique stationary distribution of $ R $ and $ B ( R , P ) $ is the zero matrix.

\end{itemize}

\end{theorem}
{\bf Proof:}
To see that $ \mu $ is an invariant measure with respect to $ R$,
note that, 
\[
\Pi_{P} + ( I - P ) D_{P} = I.
\]
Multiplying the above equation from the left by $ \mu$, yields
\begin{equation}\label{eq:erst1}
\pi_{P}^\top +  \mu^\top ( I - P) D_{P}  = \mu^\top .
\end{equation}
By simple algebra,
\begin{eqnarray}
\mu^\top & =  & \pi_{P}^\top  \sum_{ k=0}^\infty ( ( R - P )  D_P )^k  \nonumber
\\
& = & \pi_{P}^\top +   \pi_{P}^\top \sum_{ k=1}^\infty ( ( R - P )  D_P )^k   \nonumber
\\
&  = &  \pi_{P}^\top +   \pi_{P}^\top \sum_{ k=0}^\infty ( ( R - P )  D_P )^k  ( R - P  ) D_{P}  \nonumber
\\
& = & 
\pi_{P}^\top + \mu^\top ( R - P  ) D_{P}   . \label{eq:zweite2} 
\end{eqnarray}
Subtracting (\ref{eq:erst1}) from (\ref{eq:zweite2}) yields
\[
\mu^\top (  I -  R  ) D_{P}  = 0 .
\]
Existence of $ D_{P}  $ implies that $ D_{P}  = ( I - P + \Pi_{P} ) ^{-1} - \Pi_{P}  $, see (\ref{eq:ddd}).
Since $  (  I -  R  ) \Pi_{P} = 0 $, it holds that
\[
\mu^\top (  I -  R )   ( I - P + \Pi_{P} ) ^{-1} = 0 .
\]
Multiplying the above equation from the right with $  ( I - P + \Pi_{P}  ) $  yields
$ \mu  = \mu  R $, which shows that $ \mu $ is invariant to $ R$. Further, multiplying \eqref{eq:erst1} from the right with an appropriate column vector of ones, i.e., $\bar{1}$, shows
\begin{equation}
\pi_{P}^\top \bar{1} +  \mu^\top ( I - P) D_{P} \bar{1}  = \mu^\top \bar{1}  \Leftrightarrow \mu^\top \bar{1} = 1
\end{equation}
since $( I - P) D_{P} \bar{1} = ( I - \Pi_P) \bar{1}$ = 0. This shows that $\mu$ sums up to $1$.
Provided that $ \mu $ is component-wise a non-negative vector,  $ \mu $ is a stationary distribution,
which proves part (i).

For part (ii), note that by Proposition~\ref{cor:pos} it follows that $ R $ is positive recurrent. 
This together with the assumption that $ R $ is irreducible and aperiodic implies that  $ R$ is ergodic and 
\begin{equation}\label{eq:nmb}
\lim_{ n \rightarrow \infty } R^n = \Pi_R ,
\end{equation}
 where $ \Pi_R $ is a matrix with all rows equal to $ \pi_R^\top $ and $ \pi_R $ is  the unique stationary distribution of $ R $.
 Since all rows of $ \Pi_R $ are identical to $ \pi_R^\top $ and $ \mu^\top \bar{1}  = 1 $, it holds that 
 \begin{equation}\label{eq:nmb2}
 \mu^\top \Pi_R = \pi_R^\top .
 \end{equation}
We have already shown that $ \mu $  is an invariant distribution of $ R $. This together with (\ref{eq:nmb}) and (\ref{eq:nmb2}) yields
 \[
\mu^\top =  \lim_{ n \rightarrow \infty } \mu^\top   R^n = \mu^\top \Pi_R = \pi_R^\top .
 \]
 Uniqueness of the solution follows from ergodicity of $ R $ and the bias term is consequently the zeros vector, which concludes the proof.
 \hfill $ \Box $ 

\begin{remark}
Part (i) of Theorem~\ref{th:stability}  applies in case that $ R $ is a multi-chain with transient states.
In this case the stationary distribution is not unique.
This can be nicely explained via the bias term.
As the bias term depends on $ P$, it carries information on the Markov chain that 
is used in approximating  $ \pi_R $. 
Letting $ P $ tend to $ R $, the limit of $ B ( R, P ) $ typically will not tend to zero if $ R $ is a multi-chain.
This phenomenon is studied in the literature on singular perturbations, see, for example, \cite{sp1,sp2,sp3}. 

Note that uniqueness of the stationary distribution can only be established under the conditions put forward in part (ii) of  Theorem~\ref{th:stability}.
\end{remark}

The series in (\ref{eq:basicseries})  can be facilitated for deriving perturbation bounds by
\begin{eqnarray}
\pi_R^\top - \pi_ P^\top & =  & \pi_P^\top \sum_{ k=1}^\infty
( (R - P ) D_P )^k   + B ( R, P ) \label{proof} \\
& = & \pi_P^\top   ( R - P ) D_P
 \sum_{ k=0}^\infty
( ( R- P ) D_P )^k    + B ( R, P ) \nonumber \\
& = & \pi_P^\top   ( R - P ) D_P    ( I - ( R - P ) D_P )^{-1}  + B ( R, P ) 
\label{eq:SSBB}
 .
\end{eqnarray}
Following the above line of equations,
bounding $ \pi_R^\top - \pi_ P^\top $ requires bounding $     ( I - ( R - P ) D_P )^{-1} $.
We will show that the conditions put forward in the following lemma not only imply
norm bounds for  $     ( I - ( R - P ) D_P )^{-1} $ but also imply that $ B ( R , P ) $ is the zero matrix.

\begin{lemma}\label{le:basicbound}
For any matrix norm it holds with the above notation that:
\begin{itemize}

\item[(i)]
If $ || ( R - P ) D_P  || < 1 $, then
\[
||  ( I - ( R - P ) D_P )^{-1} ||  \leq \frac{1}{1 - || ( R - P ) D_P || } ,
\]

\item[(ii)] if $ || R - P  || \, ||  D_P  || < 1 $, then
\[
||  ( I - ( R - P ) D_P )^{-1} ||  \leq \frac{1}{1 - || R - P || \, || D_P ||  } ,
\]

\item[(iii)]
if $\|T\| + || R - P ||( 1 + || \pi_P ^\top || ) < 1 $, then
\[
||  ( I - ( R - P ) D_P )^{-1} ||  \leq \frac{1 - || T || }{1 - || T || - || R - P || ( 1 + ||\pi_P^\top|| )} .
\]
\end{itemize}

In addition, any 
of the conditions (i), (ii) or (iii) implies that the bias term equals the zeros vector.
\end{lemma}
{\bf Proof:}
We only provide a proof of part (iii) as the proofs of (i) and (ii) can be obtained from a similar (and simpler) line of arguments.
Using the taboo kernel representation in (\ref{eq:D}) it holds that 
\[
( R - P ) D_P  = ( R - P ) \sum_{k=0}^\infty T^k ( I - \Pi_P )  .
\]
By the condition it follows that $ || T || <1 $ and thus applying norms yields
\begin{equation}\label{eq:proof1}
|| ( R - P ) D_P ||  \leq || R - P || \frac{ 1 + || \pi_P^\top ||  }{1 - ||T|| }  .
\end{equation}
Our  condition $\|T\| + || R - P ||( 1 + || \pi_P ^\top || ) < 1 $ is equivalent to the expression on the above RHS being strictly less than 1. 
This implies that  the Neumann series  $ \sum_{k=0}^\infty (  ( R - P ) D_P  )^k  $ converges.
Consequently $ I - ( R - P ) D_P  $ is invertible with norm bounded by
\begin{eqnarray*}
||  ( I - ( R - P ) D_P )^{-1}  || & \leq & \sum_{ k=0}^\infty ||  ( R - P ) D_P  ||^k 
\\
& = & \frac{1}{1 -  ||  ( R - P ) D_P  || } .
\end{eqnarray*}
Inserting the bound in (\ref{eq:proof1}) in the expression on the above RHS concludes the proof of the statement.

For the proof of the last part of the lemma, note that $ || \pi_ R^\top  (( R - P ) D_ P )^N  || \leq  || \pi_R^\top ||  \, || ( R - P ) D_P ||^N $,
so that $ || (  R - P ) D_P ||  < 1 $ implies convergence of $ ||  \pi_R ^\top (( R - P ) D_P )^n || $ to zero as $ n $ tends to infinity. 
\hfill $ \Box $

\begin{remark} 
It is worth noting that $ || ( R - P ) D_P ||  < 1 $  typically fails in case $ R $ is a multi-chain.
Put differently, while in principle the results in the remainder of this article apply to $ R $ being a multi-chain,
we have found no example of a pair $ R , P $ with $ R $ a multi-chain and $ P $ a uni-chain such that  $ || ( R - P ) D_P ||  < 1 $.
We conjecture that  $ || ( R - P ) D_P ||  < 1 $  rules out the case that $ R $ is a multi-chain but we have not been able to prove this  so far.

\end{remark}

Note that
\[
 || ( R - P ) D_P  || \leq   ||  R - P  || \, ||  D_P  ||  \leq \frac{|| R - P || ( 1 + || \pi_P ^\top|| )}{ 1 - || T ||}
\]
implies that the bounds put forward in Lemma~\ref{le:basicbound} are increasingly limited in their applicability, while
the evaluation of the bounds becomes simpler. 
In fact, computing $ || ( R - P ) D_P  || $ is often not feasible as $ D_P $ is either
not known in closed form or is prohibitively complex in general, see  \cite{heder4,heder3,Koole}.
For the Markov chains in Example~\ref{ex:mc}, $D_P$ is known in explicit form.
For this type of problems it makes sense to apply the norm bound put forward in Lemma~\ref{le:basicbound} {\em (i)}
to (\ref{eq:SSBB}).  
More specifically, assuming $||( R - P ) D_P  || < 1$  let
\[
\Delta_{\rm DB} ( R , P ) :=  \frac{||  \pi_P^\top   ( R  - P) D_P || }{1 - ||( R - P ) D_P  ||   } ,
\]
then 
\begin{equation}\label{eq:direct}
||  \pi_R^\top - \pi_ P^\top || \leq  \Delta_{\rm DB} ( R , P )  ,
\end{equation}
which we will call the {\em direct bound} (DB).

\begin{remark}\label{rem:3}
The bound in (\ref{eq:direct}) has the following nice feature.
Let $ P $ and $ R $ be two Markov chains with $ P \not = R $ but with the same stationary distribution.
Then,   (\ref{eq:direct}) detects this and yields the correct value 0, whereas condition number type bounds yield a non-zero bound.
\end{remark}

The next bound can serve as alternative in case $D_P$ is difficult to find. It follows from replacing $(R-P)D_P$ in \eqref{eq:direct} with the taboo kernel representation and bounding the result via \eqref{eq:proof1}. Specifically, this leads to
\begin{eqnarray}
|| \pi_R^\top - \pi_ P^\top ||  & \leq   &
|| \pi_P^\top   || \,  ||  R - P  ||  \frac{ 1+ || \pi_P^\top||}{ 1 - || T || }
\frac{1  - || T||}{ 1 - || T|| - || R  - P|| (1+ || \pi_P^\top|| ) } .
\end{eqnarray}
Let
\begin{equation}
\Delta_{\rm SSB} ( R , P ) :=
|| \pi_P^\top   || \,  ||  R - P  ||  \frac{1+ || \pi_P^\top||}{ 1- || T|| - || R  - P|| ( 1+ || \pi_P^\top|| )}
 \label{eq:SSB}  ,
\end{equation}
provided that $ || T || + || R  - P|| ( 1+ || \pi_P^\top|| ) < 1 $.
Then, 
\[
|| \pi_R^\top - \pi_ P^\top || \leq \Delta_{\rm SSB} ( R , P ) 
\]
and the bound $ \Delta_{\rm SSB} ( R , P )  $ in (\ref{eq:SSB}) is called {\em Strong Stability Bound} (SSB) in the literature \cite{karta96}. 
For applications of SSB, we refer to \cite{abbassm,karta83,rabta-4,bouammor,boukir,lekadir}.
An obvious improvement of the bound in (\ref{eq:SSB}) is to replace $ || \pi_P^\top   || \,  ||  R - P  || $ by $ || \pi_P^\top   ( R - P  ) ||$; see Remark~\ref{rem:3}.

While $ P $ and $ \pi_P $ are fixed, and $ T $ offering in practice only limited flexibility, $ R $ is a
free variable of the perturbation bound.
Essentially, the direct bound and SSB only apply if $ R $ is not too far away from $ P$, i.e., if
$ || R - P || $ is small.
This is the major drawback of this type of perturbation bounds compared to condition number bounds.
To overcome this drawback, we may scale the perturbation such that the perturbation bounds do apply.
To see this,
consider the scaled model in (\ref{eq:convex}), where the
static perturbation is replaced by a scaled one, i.e., we perturb $ P $ by $ \theta ( R - P ) $ and denote the resulting transition matrix by $ P ( \theta ) $. Now, $ \theta  $ can be chosen such that the norm bounds apply to $  \theta || R - P || $.
For example, the condition on the applicability for SBB in (\ref{eq:SSB}) translates to
\[
|| T || + \theta || R  - P|| ( 1+ || \pi_P^\top|| ) < 1
\quad \text{ iff } \quad
0 \leq \theta  < \frac{ 1 - || T ||}{  || R  - P|| ( 1+ || \pi_P^\top || ) } .
\]
We call the upper bound for $ \theta $ on the RHS above the {\em domain of  SBB  with respect to $ R $}.

In the following we take an alternative route for obtaining a perturbation bound.
Starting point is (\ref{updateD}) but other than for deriving (\ref{eq:basicseries}) we now only perform the insertion operation $ K $ times,
leading to
\begin{equation}\label{eq:rty}
 \pi_{ P ( \theta )  }^\top  = \pi_{P }^\top \sum_{k=0}^K (\theta  ( R - P ) D_P)^k \,  + \pi_{P(\theta)}^\top ( \theta ( R - P) D_{ P } )^{ K+1 } .
\end{equation}
For $ K \geq 1 $, equation (\ref{eq:rty})  yields the following bound:
\begin{equation}\label{eq:seab}
\|  \pi_{ P ( \theta )  }^\top   - \pi_P^\top \|  \leq \left  \| \pi_P^\top \sum_{k=1}^K (\theta  ( R - P ) D_P)^k  \right \| \,  + \| \pi_{P(\theta)}^\top ( \theta ( R - P) D_{ P } )^{ K+1 } \| .
\end{equation}
Obviously, $ \pi_{P(\theta)}^\top $ is not known and for the actual bound we use the fact that 
\begin{eqnarray*}
  \| \pi_{P(\theta)}^\top ( \theta ( R - P) D_{ P } )^{ K+1 } \|  & \leq & 
\| \pi_{P(\theta)}^\top \|   \| ( \theta ( R - P) D_{ P } )^{ K+1 } \|  \\
& \leq & c_{||\cdot||}  \| ( \theta ( R - P) D_{ P } )^{ K+1 } \| , 
\end{eqnarray*}
where we define the norm dependent upper bound $c_{||\cdot||}$ for $\| \pi_{P(\theta)}^\top \| $ as follows
\begin{equation}\label{eq:cbound}
c_{||\cdot||} = \sup_{Q \in \mathbb{P}(S)} \| \pi_Q^\top \|,
\end{equation}
where $\mathbb{P}(S)$ represents all stochastic matrices defined on $S$. In case the $1$-norm (resp., infinity-norm) is applied to $ \pi_{P(\theta)}^\top $ we thus have 
\begin{equation}
\| \pi_{P(\theta)}^\top ( \theta ( R - P) D_{ P } )^{ K+1 } \| \leq \| ( \theta ( R - P) D_{ P } )^{ K+1 } \|.
\end{equation}
For the general $ v $-norm, a bound $c_{||\cdot||}$ can be obtained from (\ref{eq:taboo2}).

The \emph{series expansion perturbation bound of order $K$} (SEB($K$))  is now introduced by 
\begin{equation}\label{eq:SEAPertBound}
\Delta_{\rm SEB (K)} ( P ( \theta ) , P ) :=
 \left  \| \pi_P^\top \sum_{k=1}^K (\theta  ( R - P ) D_P)^k  \right \| \,  + c_{||\cdot||} \| ( \theta ( R - P) D_{ P } )^{ K+1 } \|,
\end{equation}
where $c_{||\cdot||}$ is as defined in (\ref{eq:cbound}),
and it holds that 
\[
\|  \pi_{ P ( \theta )  }^\top   - \pi_P^\top \|  \leq \Delta_{\rm{ SEB} (K)} ( P ( \theta ) , P ) ,
\]
for $ \theta \in [ 0 , 1 ] $.

\begin{remark}\label{re:ohoh}
Note that we may bound \eqref{eq:SEAPertBound} as follows
\begin{equation}\label{eq:seabNumEfficient}
\|  \pi_{ P ( \theta )  }^\top   - \pi_P^\top \|  \leq \sum_{k=1}^K \|\pi_P^\top (( R - P ) D_P)^k  \| \theta^k \,  + c_{||\cdot||} \| (( R - P) D_{ P } )^{ K+1 } \|  \theta^{K+1} ,
\end{equation}
so that the polynomial terms only have to be calculated once and can be used for evaluating the bound for different values of $\theta$. 
This is allows for fast computation and memory efficiency but, due to the additional bounding, the numerical quality of the bound decreases.
\end{remark}

From
\[
\|  ( ( R - P) D_{ P } )^{ K+1 } \| \leq \|   ( R - P) D_{ P }  \| ^{ K+1 }
\]
it follows that the series in (\ref{eq:basicseries})  converges for $P(\theta)= P + \theta ( R - P ) $ at least for $ \theta < ( \|   ( R - P) D_{ P }  \| )^{-1}  $.
Hence, for $ \theta $ sufficiently small
\begin{equation}\label{eq:series}
 \pi_{P }^\top \sum_{k=0}^K ( \theta ( R - P ) D_P)^k
\end{equation}
provides an approximation of $ \pi_{P(\theta)} $, where the error is bounded by some constant times $ \theta^{K+1}  \|  ( ( R - P) D_{ P } )^{ K+1 } \| $.
The series put forward in (\ref{eq:series}) is called  series expansion approximation  of order $ K$.
Letting $ K $ tend to infinity in (\ref{eq:series}) we obtain that
\begin{equation}\label{eq:neu3}
\pi_{ P ( \theta )}^\top =  \pi_{P }^\top \sum_{k=0}^\infty  \theta^k ( ( R - P ) D_P)^k  ,
\end{equation}
for $ 0 \leq \theta <  || ( R - P ) D_P || $.
Note that the above series expansion implies that $ \pi_{ P ( \theta )} $ tends to $ \pi_P $ as $ \theta $ tends to zero; for more details we refer to  \cite{heidergott2,heder3}.

To test the performance of the different bounds in the scaled perturbation setting (i.e., \eqref{eq:convex}) we will investigate the relative error of the perturbation bounds. Clearly, a better bound results in a smaller relative error. Consider a condition number bound for $ || \pi_{ P ( \theta )}^\top - \pi_P^\top || $.
The following reasoning only uses the basic definition of a CNB in (\ref{opl}) so that the arguments apply to the
condition number bounds for finite chains discussed in Section~\ref{sec:norms2} and the CNB in Remark~\ref{re:1} as well.
Generally speaking, let $ \Delta_{\rm CNB} ( P(\theta) , P ) =  \theta \kappa || R - P || $ denote a condition number bound for $ \| \pi^\top_{P(\theta)} - \pi_P^\top \|$.  
Following (\ref{relerr}),
the relative error inferred by using $ \theta \kappa || R - P ||  $ rather than $ || \pi_{ P ( \theta )}^\top - \pi_P^\top ||  $
is given by
\begin{eqnarray}\label{eq:err}
\eta_{\rm{CNB}}(\theta)  & :=& 
\frac{ \Delta_{\rm CNB} ( P ( \theta) , P )   - || \pi_{ P ( \theta )}^\top - \pi_P^\top ||  }{ || \pi_{ P ( \theta )}^\top - \pi_P^\top ||  } 
\nonumber \\
& = & 
\frac{ \theta  \kappa|| R - P || - || \pi_{ P ( \theta )}^\top - \pi_P^\top ||}{ || \pi_{ P ( \theta )}^\top - \pi_P^\top || }
\nonumber
\\
&  = & \frac{ \theta \kappa|| R - P || }{ || \pi_{ P ( \theta )}^\top - \pi_P^\top || } - 1.
\end{eqnarray}
Note that this relative error is by definition $\geq 0$. 
In the same vein, when we replace $ {\bf \Delta } ( R , P ) $ in (\ref{relerr}) by the bounds $  \Delta_{\rm{SSB}}(P ( \theta) , P)$,
$  \Delta_{\rm{DB}}(P ( \theta) , P) $ and $  \Delta_{\rm{SEB}(K)}(P ( \theta) , P) $, respectively, we obtain the corresponding absolute relative error expressions  denoted by $  \eta_{\rm{SSB}}( \theta)$,
$  \eta_{\rm{DB}}(\theta) $ and $  \eta_{\rm{SEB}(K)}( \theta) $.

The following theorem analyses the relative error of the discussed bounds. It shows that in general the relative error of a condition number bound and SSB  converges for $\theta \downarrow 0$ to a finite non-zero value, while the SEB($K$)-based bounds have the desirable property that the relative error vanishes. 
Moreover, the rate of convergence of the relative error of SEB($K$) can be explicitly computed.

\begin{theorem}[Relative Errors]\label{th:relerror}
Let  $\| \pi_{P(\theta)}^\top - \pi_P^\top \| > 0$, for all $\theta \in (0,1]$. 

\begin{itemize}

\item[(i)] The relative error of  the condition number bound (CNB)  is given by
\[
\eta_{\rm{CNB}}(\theta) = \frac{ || R - P || \kappa}{ || \pi_{P(\theta)}^\top  ( R - P ) D_P  || } - 1,
\]
and it holds that  
\[
\lim_{\theta \downarrow 0} \eta_{\rm{CNB}}(\theta) = \frac{ || R - P || \kappa}{ || \pi_P^\top  ( R - P ) D_P  || } - 1 \geq 0 ,
\]
where equality is only reached in the special case when $|| R - P || \kappa$ equals $|| \pi_P^\top  ( R - P ) D_P  ||$.

\end{itemize}

\noindent

\begin{itemize}
\item[(ii)] Provided that $ || T || +  \theta || R  - P|| ( 1+ || \pi_P^\top|| ) < 1 $, the relative error of the strong stability bound (SSB) is given by
\[
\eta_{\rm{SSB}}(\theta) = 
\frac{ ||  R - P  || \, || \pi_P^\top   ||  (1+ || \pi_P^\top ||) }{ || \pi_{P(\theta)}^\top (R - P) D_P  ||  (1 - || T || - \theta || R  - P|| ( 1+ || \pi_P^\top || ))} -1,
\]
and it holds that  
\[
\lim_{\theta \downarrow 0} \eta_{\rm{SSB}}(\theta) =  \frac{ ||  R - P  || \, || \pi_P^\top   ||  (1+ || \pi_P^\top ||) }{ || \pi_P^\top (R - P) D_P  ||  (1 - || T ||)} -1 \geq 0 ,
\]
where equality is only reached in the special case when the nominator equals the denominator in the fraction.

\item[(iii)] Provided that  $\theta\|(R-P)D_P\| < 1$, the relative error of the direct bound (DB) is given by
\[
\eta_{\rm{DB}}(\theta) = \frac{\frac{\|\pi_P^\top (R-P)D_P\|}{1-\theta\|(R-P)D_P\|}}{ || \pi_{P(\theta)}^\top (R - P) D_P  ||} - 1 ,
\]
and it holds that  $\lim_{\theta \downarrow 0} \eta_{\rm{DB}}(\theta) =  0$.

\item[(iv)] Provided that $\theta\|(R-P)D_P\| < 1$, the relative error of the  series expansion bound of order $ K \geq 1$ (i.e., {\rm{SEB}}($K$)) is given by 
\[
\eta_{\rm{SEB}(K)}(\theta) = \frac{  2 c_{||\cdot||}  \| ( ( R - P) D_{ P } )^{ K+1 } \| \theta^K  }
{ || \pi_{P(\theta)}^\top ( R - P ) D_P|| } ,
\]
and it holds that 
$ \eta_{\text{SEB($K$)}}(\theta) $ is of order $ O ( \theta^{K} )$.

\end{itemize}
\end{theorem}
{\bf Proof:}
All  relative error expressions follow by simply inserting the different bounds and using the result that
\begin{equation}\label{eq:partProofRelErrors}
\pi_{ P ( \theta )}^\top -  \pi_{P }^\top =  \theta \pi_{P(\theta)}^\top ( R - P ) D_P
\end{equation}
in the denominator is of order $ O ( \theta )$.
Indeed, using the fact that $ P ( \theta ) $ is irreducible and aperiodic for $ \theta < 1 $  it follows from  (\ref{eq:basicseries}) together with Theorem~\ref{th:stability} that 
\begin{equation}\label{eq:split}
 \pi_{P(\theta)}^\top ( R - P ) D_P   = \frac{1}{\theta} (   \pi_{P(\theta)}^\top -  \pi_{P}^\top)
  =  \pi_P^\top ( R -P ) D_P + \pi_P^\top \sum_{k=2}^\infty \theta^{k-1}  ( ( R -P ) D_P )^k  ,
\end{equation}
which shows that $  \pi_{P(\theta)}^\top ( R - P ) D_P  $ can be written as  power series with leading term 
$ \pi_P^\top  ( R -P ) D_P $, and thus 
 implies that $ \theta \|  \pi_{P(\theta)}^\top ( R - P ) D_P \|  $ is of order $ O ( \theta ) $.

We now turn to perturbation bounds.
For CNB it holds that
\[
\eta_{\text{CNB}}(\theta) = \frac{ \theta || R - P || \kappa }{ || \pi_{ P ( \theta ) }^\top - \pi_P^\top  || } - 1
=
\frac{ || R - P || \kappa}{ || \pi_{P(\theta)}^\top  ( R - P ) D_P ) || } - 1
\]
and the limit result then follows from (\ref{eq:split}), where the second equality is obtained by  (\ref{eq:partProofRelErrors}).

The proof of the statements for SSB, DB and SEB($K$) follow from the same line of argument and we will 
in the following only present the proof for the most challenging of these cases which is the relative error of the $K$-th order SEB.
Following (\ref{eq:SEAPertBound})  we can write
\begin{equation}\label{eq:relErrorSE}
\eta_{\text{SEB($K$)}}(\theta) = \frac{   \overbrace{ \left \| \pi_P^\top \sum_{k=1}^K (\theta  ( R - P ) D_P)^k  \right \|}^{=:H} \,  +  c_{||\cdot||} \| ( \theta ( R - P) D_{ P } )^{ K+1 } \|  }
{ || \pi_{ P ( \theta ) }^\top - \pi_P^\top  ||} - 1 .
\end{equation}
For $ H $ it holds that 
\[
H = \left \| \pi_P^\top \sum_{k=0}^{K-1} (\theta  ( R - P ) D_P)^k \theta  ( R - P ) D_P \right \|.
\]
After some algebra,
\[
H = \left \| \pi_P^\top \sum_{k=0}^\infty (\theta  ( R - P ) D_P)^k \left[ I - (\theta  ( R - P ) D_P)^K \right] \theta ( R - P ) D_P \right \|
\]
and using condition $\|\theta(R-P)D_P\| < 1$ together with (\ref{eq:neu3})  
we arrive at 
\[
H = \left \| \pi_{P(\theta)}^\top \left[ I - (\theta  ( R - P ) D_P)^K \right] \theta ( R - P ) D_P \right \| ,
\]
which can be straightforwardly bounded by
\[
H \leq \| \pi_{P(\theta)}^\top \theta ( R - P ) D_P \| + c_{||\cdot||} \| (\theta  ( R - P ) D_P)^{K+1}\|.
\]
Inserting the above bound for $H$ into \eqref{eq:relErrorSE} yields for the relative error
\begin{equation}\label{eq:boundRelErrorSEBK}
\eta_{\text{SEB($K$)}}(\theta) \leq \: \theta^{K+1} \, \frac{ 2 c_{||\cdot||} \| (( R - P ) D_P)^{K+1}\|  }
{ || \pi_{ P ( \theta ) }^\top - \pi_P^\top  || }.
\end{equation}
The limit results now follows from the fact that $ || \pi_{ P ( \theta ) }^\top - \pi_P^\top  || $  is of order $ O ( \theta ) $.
\hfill $ \Box $

\begin{figure}[ht]
\begin{picture}(400,250)
\includegraphics[scale=0.8]{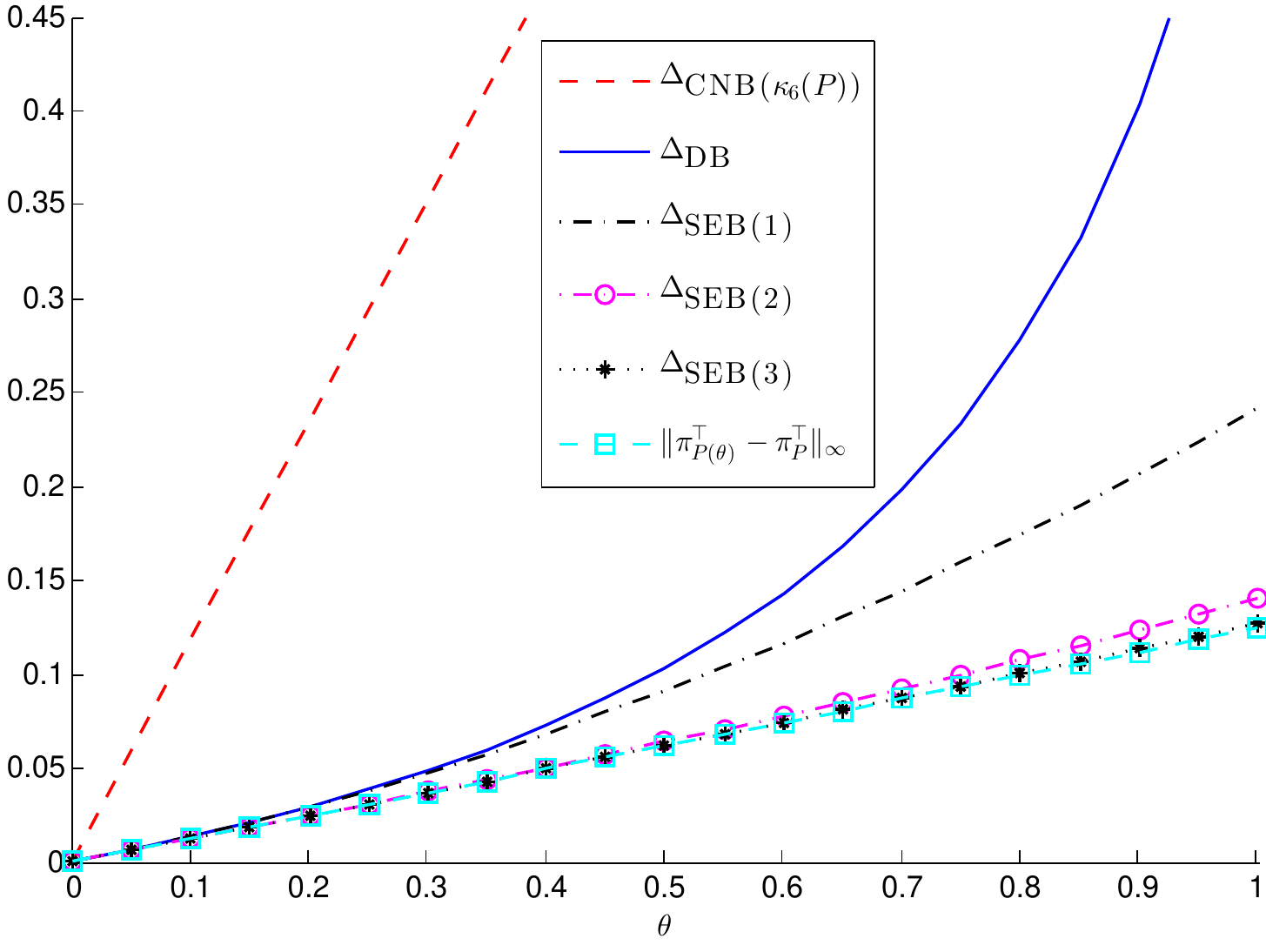}
\end{picture}
\caption{Perturbation bounds for $\|\pi_{P(\theta)}^\top-\pi_P^\top\|_\infty$ with $\theta \in (0,1]$, where $P(\theta)=(1-\theta)P+\theta R$ for randomly generated $P$ and $R$ consisting of $40$ states.}
\label{fig:pertBoundsSimpleExample}
\end{figure}

\begin{figure}[ht]
\begin{picture}(400,250)
\includegraphics[scale=0.8]{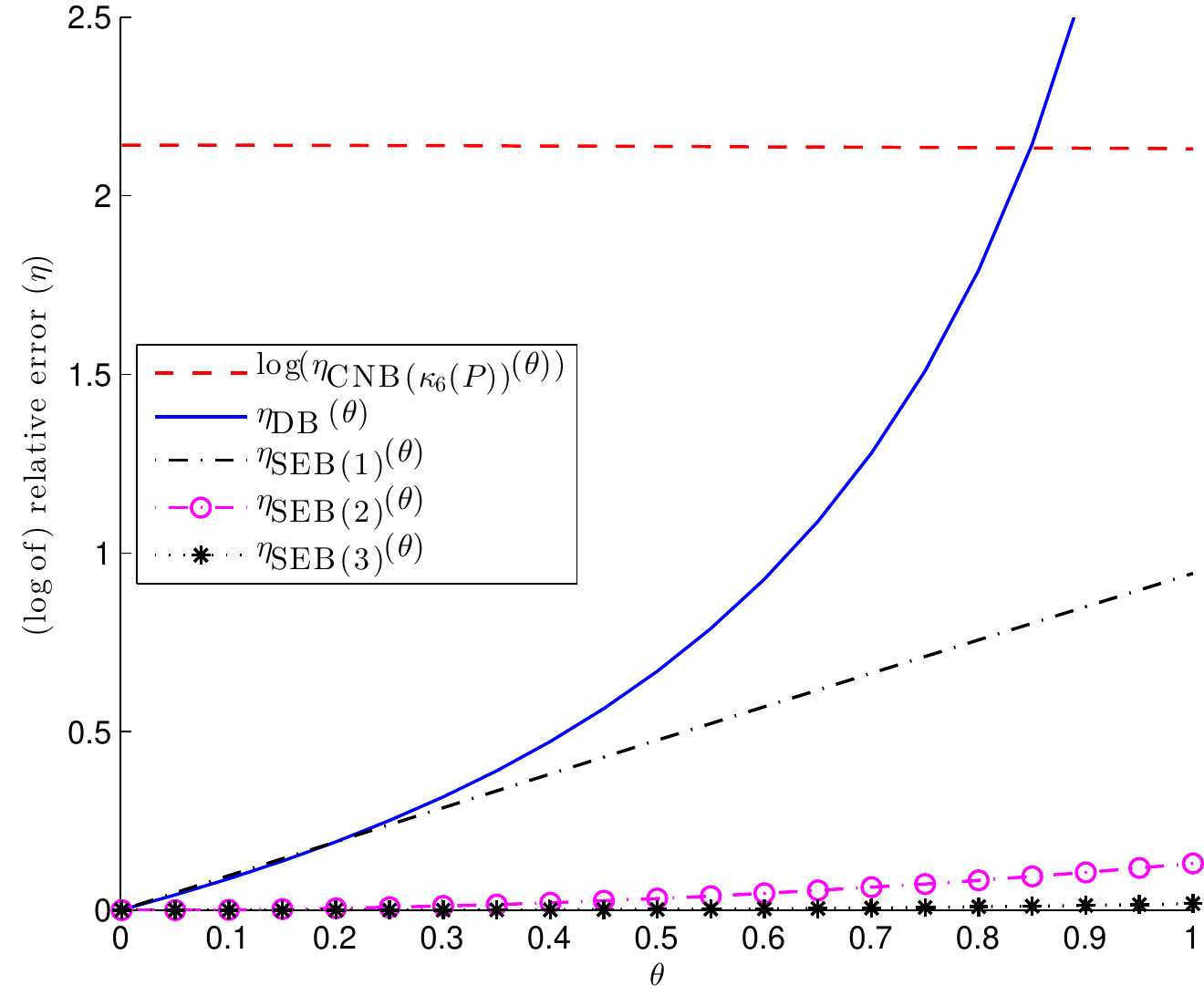}
\end{picture}
\caption{Relative errors of the perturbation bounds for $\|\pi_{P(\theta)}^\top -\pi_P^\top \|_\infty$ with $\theta \in (0,1]$, where $P(\theta)=(1-\theta)P+\theta R$ for randomly generated $P$ and $R$ consisting of $40$ states.}
\label{fig:pertBoundsSimpleExampleRelErrors}
\end{figure}

For an illustration 
of Theorem~\ref{th:relerror} we generated two random transition matrices $P$ and $R$ with $40$ states. The random generation is done by drawing random numbers from $(0,1)$ and normalizing the rows so that they sum up to 1. 
Then we considered in case of the $\infty$-norm all perturbation bounds from Theorem~\ref{th:relerror} on the interval $\theta \in (0,1]$ together with the true perturbation effect $\| \pi_{ P ( \theta )}^\top - \pi_P^\top \|_\infty$ (the true effect was calculated numerically). 
The results can be found in Figure~\ref{fig:pertBoundsSimpleExample}. 
Figure~\ref{fig:pertBoundsSimpleExample} shows that in this experiment all bounds, except for CNB, are similar in performance on the interval $\theta \in [0,0.1]$. 
For $\theta > 0.1$ SEB of order $K=3$ performs best.
DB performs similar to SEB($1$) on the interval $\theta \in (0,0.3]$ but for $\theta > 0.3$ SEB($1$) outperforms DB. 
This simple example illustrates that in a scaled perturbation setting CNB is apparently too general to be competitive compared to the other bounds. 
The differences become more apparent if we look at the relative errors for the different bounds plotted in Figure~\ref{fig:pertBoundsSimpleExampleRelErrors}. 
The results for SSB are not plotted because the condition in part (ii) of Lemma~\ref{le:basicbound} is not met.

\begin{remark}\label{re:6}
Provided that $\theta_0$ exists such that $\theta_0 || ( R - P ) D_P  || < 1$,
then 
\[
\eta_{\text{SEB($K$)}}(\theta) = O(\theta_0^{K}) ,
\]
for $ 0 \leq \theta \leq \theta_0 $.
\end{remark}

\begin{remark}
The result put forward in Theorem~\ref{th:relerror} seems to contradict the fact that for finite Markov chains it
holds that 
\begin{equation}\label{eq:relrel}
\left | \frac{ (\pi_R )_i - (\pi_P )_i }{ (\pi_R  )_i  } \right | \leq 2 \eta n + O ( \eta^2 ) , \: i \in S = \{ 0 , \ldots ,n-1 \} ,
\end{equation}
where $ \eta $ is bounded by $ || R - P || $  and $n$ denotes the size of the state-space,
which indicates that the relative element-wise error in using $ \pi_P $ as a substitute for $ \pi_R $ tends to zero 
as $ P $ approaches $ R $, see \cite{n1,n2,n3} for details.
Note that above equation is equivalent to 
\[
\left |  (\pi_R )_i - (\pi_P )_i \right | \leq   (\pi_R  )_i  \left (  2 \eta n + O ( \eta^2 ) \right ) , \: i \in S = \{ 0 , \ldots ,n-1 \} ,
\]
and reads in norm-version, using, for example, the $\infty$-norm (or 1-norm), 
\[
|| \pi_R^\top - \pi_P^\top ||_{ \bf 1} \leq  2 \eta n + O ( \eta^2 )  ,
\]
see Remark~\ref{re:ohoh}. 
Hence, the element-wise  relative error result in (\ref{eq:relrel}) is a statement about continuity of finite Markov chains and does not imply that 
the relative error in predicting the true norm distance between $ \pi_R $ and $ \pi_P $ by a CNB becomes small; for details
compare the  definition of the relative error in (\ref{relerr}) and (\ref{eq:err}), respectively, with that in (\ref{eq:relrel}).

\end{remark}

We conclude this section by presenting an interesting result for stability theory.

\begin{corollary}\label{cor:suff}
Consider the model  $ P ( \theta ) = ( 1- \theta ) P + \theta R $, $ \theta \in [ 0 ,1 ) $, with $ P $ aperiodic, irreducible and positive recurrent.
If 
\[
\theta < \frac{1 - || _i \! P  ||}{ || R - P ||} ,
\]
then $ P ( \theta) $ has a unique stationary distribution.
\end{corollary}
{\bf Proof:}
Note that  $ P ( \theta ) $ is aperiodic and irreducible  for $ \theta \in [ 0 , 1 )$.
It remains to be shown that $ P ( \theta ) $ is positive recurrent.
By computation,
\begin{eqnarray*}
|| _i ( P ( \theta ) )  || & = & 
 || _i (  ( 1- \theta ) P + \theta R ) || \\
  & \leq  & || _i \! P  + \theta ( R - P ) || 
 \\
 &  \leq &  || _i \! P  || + \theta || R - P || .
 \end{eqnarray*}
Hence, provided that $ \theta $  satisfies $  || _i \! P  || + \theta || R - P || < 1 $, it follows $ || _i ( P ( \theta ) ) || < 1  $ and 
by  Proposition~\ref{cor:pos} we conclude that $ P ( \theta ) $ is positive recurrent.
Solving $ \theta $ out of $  || _i \! P  || + \theta || R - P || < 1 $ concludes the proof.
\hfill $ \Box $ 

\begin{remark}\label{re:ssb}
Note that from Corollary~\ref{cor:suff} it follows that  if condition (ii) in Theorem~\ref{th:relerror} for the SSB with $T = {}_i P $, for some $i \in S$,  is satisfied, 
then $ P ( \theta ) $ is stable, i.e., has a unique stationary distribution.
\end{remark}
Kartashov established in \cite{karta96} a result similar to Theorem~\ref{th:stability}.
It is worth noting that Kartashov didn't provide a lower bound for the region of stability as detailed in 
 Corollary~\ref{cor:suff} together with Remark~\ref{re:ssb}.

\section{Explicit Perturbation Bounds for the Two-State Chain (Finite State Space)}\label{sec:leadex}
In this section we explicitly compute the bounds put forward in Theorem~\ref{th:relerror} for the 
two-state chain from Example~\ref{ex:mc}. The following convex combination is considered
\[
P(\theta) = (1 - \theta) \underbrace{\begin{pmatrix}
      1 - p & p \\
      q & 1 - q \\
    \end{pmatrix}}_{=P^s} + \theta \underbrace{\begin{pmatrix}
      1 - \widetilde{p} & \widetilde{p} \\
      \widetilde{q} & 1 - \widetilde{q} \\
    \end{pmatrix}}_{:=\widetilde{P}^s} .
\]
We are interested perturbing $P(0)$ by choosing $\theta > 0$. Note that for the difference in Markov transition matrices it holds
\[
P(\theta) - P(0) = \theta( \widetilde{P}^s  - P^s ) = \theta \begin{pmatrix}
                      p - \widetilde{p} & \widetilde{p} - p\\
                      \widetilde{q} - q & q - \widetilde{q} \\
                    \end{pmatrix}.
\]
which gives 
\[
|| P(\theta) - P(0)  ||_{v} =  \theta (1 + \alpha ) \max \left\{ | p - \widetilde{p}|, \frac{ 1 }{ \alpha }| q - \widetilde{q}| \right\} .
\]

In the following the explicit perturbation bounds are presented for the $v$-norm. Using \eqref{updateD} in the calculation for CNB we get 
\begin{equation*}
|| \pi_{P(\theta)}^\top - \pi_{P^s}^\top ||_{v} \leq  || \pi_{P(\theta)}^\top ||_{v} ||P(\theta) - P^s||_{v} ||D_{P^s}||_{v}.
\end{equation*}
It holds that (see also Example~\ref{ex:mc})
\begin{equation*}
|| \pi_{P(\theta)}^\top ||_{v} \leq \alpha \qquad \mbox{ and } \qquad ||D_{P^s}||_v = \frac{1+\alpha}{(p+q)^2}\max \left \{ p , \frac{q}{\alpha} \right \} 
\end{equation*}
so that we obtain for CNB
\begin{equation*}
\Delta_{\rm {CNB}} (P(\theta) , P^s ) = 
\theta \left(\frac{1+\alpha}{p+q}\right)^2 \max \left\{ \alpha  |p-\widetilde{p}| , |q-\widetilde{q}| \right\} \max\left\{ p , \frac{q}{\alpha} \right\}.
\end{equation*}
In the general framework of CNB given in \eqref{opl} with \eqref{Delta} it holds that $\kappa = \frac{1+\alpha}{(p+q)^2}\max\{\alpha p , q\}$ for this example. 

For the SSB we compute
\[
  || \pi_{P^s}^\top ||_{v}=\frac{q + p\alpha}{p + q}.
\]
Next, the individual terms in  (\ref{eq:SSB}) have to be computed.
Here, we make use of the taboo kernel bound as provided in Example~\ref{ex:T}, where the taboo kernel may be obtained by removing one of the columns where the choice of column depends on the value of $ p $ and $ q $, and we arrive at 
\[
|| T||_v \leq \min \{  \max\{ \alpha p, 1 - q \} , \max \{   1-p   , q \} \} .
\]
Note that a similar analysis can be carried out when considering removing rows of $ P^s $.
SSB can only be provided for small perturbations, i.e., small values of $\theta$. More specifically, provided that 
\[
\theta < \frac{ 1 -  \min \{ \max\{ \alpha p, 1 - q \} , \max \{  \alpha ( 1-p )  , q \} }{ \left( 1 + \frac{q+p \alpha}{p+q} \right) (1 + \alpha ) \max \left\{ | p - \widetilde{p}|, \frac{ 1 }{ \alpha }| q - \widetilde{q}| \right\} } ,
\]
the SSB bound for $ || \pi_{P(\theta)}^\top - \pi_{P^s}^\top ||_{v}$  is given  by
\[
\Delta_{\rm {SBB}} (P(\theta) , P^s ) = \frac{ \left( \frac{q+p \alpha}{p+q} \right) \left( 1 + \frac{q+p \alpha}{p+q} \right) \theta (1 + \alpha ) \max \left\{ | p - \widetilde{p}|, \frac{ 1 }{ \alpha }| q - \widetilde{q}| \right\}}
{1 -      \min \{  \max\{ \alpha p, 1 - q \} , \max \{   1-p   , q \} \}            - \left( 1 + \frac{q+p \alpha}{p+q} \right) \theta (1 + \alpha ) \max \left\{ | p - \widetilde{p}|, \frac{ 1 }{ \alpha }| q - \widetilde{q}| \right\}   }.
\]
For example, letting $ \alpha = 1 $, which is possible, see Lemma~\ref{le:taboo}, yields the simplified expression
\[
\Delta_{\rm {SBB}} (P(\theta) , P^s ) = \frac{4  \theta  \max \left\{ | p - \widetilde{p}|, | q - \widetilde{q}| \right\}}{1 -
 \min \left\{ \max\{ p, 1 - q \} , \max \{   1-p   , q \}  \right\}    - 4  \theta
 \max \left\{ | p - \widetilde{p}|, | q - \widetilde{q}| \right\}
 }
\]
for SSB.
By inspection of above, it is obvious that SSB behaves poorly for $ p $ and $ q $ close to one or close to zero as in this case
the norm of the taboo kernel approaches one.

Calculations show that DB leads to
\begin{equation*}
\Delta_{\rm {DB}} (P(\theta) , P^s ) = \frac{\theta|p\widetilde{q}-\widetilde{p}q|(1+\alpha)}{(p+q) \left( p + 1 - \theta(1+\alpha) \max \{ |p-\widetilde{p}| , \frac{|q-\widetilde{q}|}{\alpha} \} \right) }
\end{equation*}
under the assumption that
\begin{equation*}
\theta < \frac{p+1}{(1+\alpha)\max \{ |p-\widetilde{p}| , \frac{|q-\widetilde{q}|}{\alpha} \} }.
\end{equation*}

For SEB($K$) with $K=0$ it holds
\begin{equation*}
\Delta_{\rm {SEB} (0) } (P(\theta) , P^s ) = \frac{\theta(1+\alpha)}{p+q} \max \{ \alpha |p-\widetilde{p}|, |q-\widetilde{q}| \}
\end{equation*}
of which the construction is similar to CNB but with the difference that CNB requires an additional bounding on $||(P(\theta)-P^s)D_{P^s}||_v$ to obtain $||(P(\theta)-P^s)||_v ||D_{P^s}||_v$, which stems from the fact that $||(P(\theta)-P^s)D_{P^s}||_v \leq ||(P(\theta)-P^s)||_v ||D_{P^s}||_v$.
More specifically, $ \Delta_{\rm {CNB}} (P(\theta) , P^s ) $ is by factor 
\begin{equation*}
\frac{\Delta_{\rm {CNB}} (P(\theta) , P^s ) }{\Delta_{\rm {SEB}(0)} (P(\theta) , P^s ) } = \frac{1+\alpha}{p+q} \max\left\{p,\frac{q}{\alpha}\right\}  \geq 1 
\end{equation*}
larger than $ \Delta_{\rm {SEB}(0)} (P(\theta) , P^s ) $. 
In case $\alpha = 1$ this factor is $2\max\{p,q\}/(p+q)$, which is greater than 1 for $p \neq q$. When $\alpha$ is chosen to be $>>1$ this factor grows linearly in $\alpha$. This illustrates that, although being more general, CNB loses on quality in contrast to SEB$(0)$ since it does not utilize the contraction property of $(P(\theta)-P^s)D_{P^s}$.

After similar calculations it can be shown that SEB($K$) with $K=1$ results in
\begin{equation*}
\Delta_{\rm {SEB} (1) } (P(\theta) , P^s )  = \frac{\theta(1+\alpha)}{(p+q)^2} \left( |p\widetilde{q} - \widetilde{p} q| + \theta | p - \widetilde{p} + q - \widetilde{q} | \max \{ \alpha |p-\widetilde{p}|, |q-\widetilde{q}| \} \right).
\end{equation*}

\section{An Elaborate Perturbation Analysis of a Queueing System}\label{sec:XXX}

To illustrate the application of perturbation bounds in a setting where the deviation matrix is not available in a closed-form, we discuss in this section
 the M/G/1 queue with breakdowns.
In addition, we consider the finite version of the queue, i.e., the M/G/1/N queue with breakdowns
and we illustrate  SEB($K$).
The breakdown model will have the special feature that we perturb the system with no breakdowns
by an unstable chain modeling a pure birth process.

The basic model of the M/G/1 queue with breakdowns is introduced in Section~\ref{sec:A1} and
in Section~\ref{sec:A3} a discussion of the literature is provided.
The perturbation bounds for both models are presented in Section~\ref{sec:A5} and Section~\ref{sec:A6}, respectively.

\subsection{The Basic Model}\label{sec:A1}

Consider a single server queue.
Customers arrive at the queue according to a Poisson-$\lambda$-arrival process. Service times are identically distributed with mean $ 1 / \mu $ and we denote the service time distribution by $ \mathcal{S}(x) $. 
Throughout this section we assume that $ \lambda / \mu < 1$.
At the beginning of each service, there is a probability $ \theta $ that the server breaks down (and
the customer is send back to the queue) and enters a repair state, the length of which is exponentially distributed with rate $ r $ and which is independent of everything else, and with probability $ ( 1 - \theta ) $ the server is operational and serves the customer (if any, according to FCFS). The only points in time where a possible server breakdown can occur is right at the beginning of a service. This system is modeled by the jump chain embedded at service completions and completions of a repair, and it has state space $S=\{0,1,\dots\}$. The transition probabilities from $i \in S$ to $j \in S$, denoted as $P_\theta(i,j)$, are given as follows:

\smallskip

For $i = 0$, the process jumps to $ j \geq 0 $ if
a customer arrives and the server is operational and during the service of this customer there are $j$ additional arrivals.
This probability is given by
\[
( 1 - \theta ) \int_{0}^{\infty} e^{ -\lambda x } \frac{( \lambda x )^{j}}{j!} \, d  \mathcal{S} ( x ).
\]
Alternatively,
a customer arrives at the empty queue and  the server breaks down at service initiation and during the repair time of the server there are $ j-1 $ additional arrivals, so that at the end of the repair time there are in total $ j$ customers at the server.
This probability is given by
\[
\theta  \int_{0}^{\infty} e^{ -\lambda x } \frac{( \lambda x )^{j-1}}{(j-1)!} \, r e^{ -r x } d  x = \theta \frac{r}{\lambda + r} \left( \frac{\lambda}{\lambda + r} \right)^{j-1},
\]
for $ j \geq 1 $ and zero for $ j=0$, where we make use of the convention that $ 0! = 1$.
Combining these results, for $ i =0  $, we arrive at
\[
P_\theta (  0 ,  j  ) = ( 1 - \theta ) \int_{0}^{\infty} e^{ -\lambda x } \frac{( \lambda x )^{j}}{j!} \, d \mathcal{S} ( x ) + \theta \frac{r}{\lambda + r} \left( \frac{\lambda}{\lambda + r} \right)^{j-1} 1_{ j \geq 1} .
\]

For $ i \geq 1 $, the process jumps to state $ j \geq i-1 $ if
the server remains operationally, so that service of the subsequent customer in the queue may begin, and
during the service of this customer there are $j - i + 1 \geq 0$ additional arrivals.
This probability is given by
\[
( 1 - \theta ) \int_{0}^{\infty} e^{ -\lambda x } \frac{(\lambda x)^{j - i + 1}}{(j - i +1)!} \, d \mathcal{S} ( x ) .
\]
Alternatively, there is a server breakdown and during the exponential repair time there are $ j - i  \geq 0 $ arrivals from the outside. This probability is given by
\[
\theta \frac{r}{\lambda + r} \left( \frac{\lambda}{\lambda + r} \right)^{j - i }.
\]
Combining these results, we arrive at
\begin{equation}\label{eq:Pcentral}
P_\theta ( i , j ) = ( 1 - \theta ) \int_{0}^{\infty} e^{ -\lambda x } \frac{(\lambda x)^{j - i + 1}}{(j - i +1)!} \, d \mathcal{S} ( x ) + \theta \frac{r}{\lambda + r} \left( \frac{\lambda}{\lambda + r} \right)^{j-i}  1_{ j  \geq i },
\end{equation}
for $ 1 \leq i $ and $ i  -1\leq j $. All other entries of $ P_\theta  $ are set to zero.

Observe that for $ \theta = 1 $, $ P_1 $ models a pure birth process and the queue is not stable, whereas $ P_0 $ models a stable M/G/1 queue with no breakdowns. 
The kernel $ P_\theta $ is given through the convex combination $ \theta P_1 + ( 1- \theta ) P_0 $ of the two kernels.

\subsection{Discussion of Literature}\label{sec:A3}

Since the pioneering work of Thiruvengadam \cite{thir} and
Avi-Itzhak and Naor \cite{avi}, there has been a considerable
interest in the study of queues with server breakdowns, see for
example \cite{cao1982,lishi,wangcao} and references therein. However,
the majority of results is expressed in terms of  systems of
equations the solution of which is rather challenging, or have
solutions which are not easily interpretable in practice. For
instance, Baccelli and Znati \cite{bacz} provide the generating
function of the number of customers in the $M/G/1$ system with
dependent breakdowns. Also, results are given in terms of the
inverse of Laplace transforms, see, e.g., \cite{bacz}, which
require numerical inversion for solving a given system. To overcome
these difficulties, approximation methods are used where  the
complex (real) system is replaced by one which is ``close'' to it
in some sense but which has a simpler structure (resp.,
components) and for which analytical results are available. 

\subsection{The Infinite Capacity M/G/1 Queue with Breakdowns  (Denumerable State Space)}\label{sec:A5}
In this section the M/G/1 queue with breakdowns is considered. Note that SSB is the only bound applicable as the size of the state-space is infinite and the deviation matrix is not known in explicit form.
As next we provide auxiliary results for obtaining the overall SSB. 
Recall that $ P_0 $ is the transition kernel of the embedded jump chain of an M/G/1 queue and we consider the taboo kernel $ T = _0 \! (P_0 )$, i.e., 
we remove the first column of $ P_0$.

For the taboo kernel $T$ it holds that 
\begin{eqnarray}
  \|  T \|_{\upsilon}
   & = & \sup\limits_{i \geq 0} \frac{ 1 }{ \alpha^{i} } \sum\limits_{j \geq 1} \alpha^{j} \left| \int_{0}^{\infty} e^{-\lambda x} \frac{(\lambda x)^{j - i +1 }}{ ( j - i + 1 )! } d\mathcal{S}(x)\right| 1_{ j - i + 1 \geq 0 } \nonumber
\\
   & = & \sup\limits_{i \geq 0}   \frac{ 1 }{ \alpha^{i} } \sum\limits_{j \geq 1} \alpha^{j} \int_{0}^{\infty} e^{-\lambda x} \frac{(\lambda x)^{j - i +1 }}{ ( j - i + 1 )! } d\mathcal{S}(x) 1_{ j \geq i-1 } \nonumber
\\
   & = & \sup\limits_{i \geq 0} \frac{ 1 }{ \alpha^{i} } \sum\limits_{j \geq \max(i-1,1)} \alpha^{j }  \int_{0}^{\infty} e^{-\lambda x} \frac{(\lambda x)^{ j }}{ j! } d\mathcal{S}(x)
 \nonumber
\end{eqnarray}
For $ i = 0,1 $,
\begin{eqnarray}
\sup\limits_{0 \leq i \leq 1} \frac{ 1 }{ \alpha^{i} } \sum\limits_{j \geq \max(i-1,1)} \alpha^{j } \left| \int_{0}^{\infty} e^{-\lambda x} \frac{(\lambda x)^{ j }}{ j! } d\mathcal{S}(x) \right|
&  = &   \sum\limits_{j \geq 1}  \alpha^{j } \int_{0}^{\infty} e^{-\lambda x} \frac{(\lambda x)^{ j }}{ j! } d\mathcal{S}(x)
 \nonumber \\
&  = &   \sum\limits_{j \geq 1}   \int_{0}^{\infty} e^{-\lambda x} \frac{(\lambda \alpha x)^{ j }}{ j! } d\mathcal{S}(x)
 \nonumber \\
&  = &     \int_{0}^{\infty} e^{-\lambda x} \sum\limits_{j \geq 1}  \frac{(\lambda \alpha x)^{ j }}{ j! } d\mathcal{S}(x)
 \nonumber \\
&  = &     \int_{0}^{\infty} e^{-\lambda x} ( e^{ \lambda \alpha x}  -1 )  d\mathcal{S}(x)
 \nonumber \\
& = &
\int_{0}^{\infty} e^{- \lambda (1- \alpha ) x} d\mathcal{S}(x) -     \int_{0}^{\infty} e^{-  \lambda  x} d\mathcal{S}(x) ,
\nonumber
\end{eqnarray}
and for $ i > 1 $ 
\begin{eqnarray}
&  & \sup\limits_{i \geq 2} \frac{ 1 }{ \alpha^{i} } \sum\limits_{j \geq \max(i-1,1)} \alpha^{j  - 1} \left| \int_{0}^{\infty} e^{-\lambda x} \frac{(\lambda x)^{ j }}{ j! } d\mathcal{S}(x) \right| \nonumber
\\
&  & \qquad =
\sup\limits_{i \geq 2} \frac{ 1 }{ \alpha^{i} } \sum\limits_{j \geq i-1} \alpha^{j  - 1} \int_{0}^{\infty} e^{-\lambda x} \frac{(\lambda x)^{ j }}{ j! } d\mathcal{S}(x)
 \nonumber \\
   &  & \qquad  =  \frac{ 1 }{ \alpha^3 } \int_{0}^{\infty} e^{-\lambda x} \sum\limits_{j \geq 0} \frac{(\lambda \alpha x)^{ j }}{ j! } d\mathcal{S}(x)
 -   \frac{ 1 }{ \alpha^3 } \int_{0}^{\infty} e^{-\lambda x}  d\mathcal{S}(x)
 \nonumber \\
   &  & \qquad = \frac{1}{\alpha^3 }
\left (  \int_{0}^{\infty} e^{ -\lambda (1 - \alpha ) x} d\mathcal{S}(x)
-   \int_{0}^{\infty} e^{-\lambda x}  d\mathcal{S}(x) \right )   \nonumber  .
\end{eqnarray}
Denoting by  $\mathcal{S}^{*}(z)$ the Laplace-Stieltjes transform of $\mathcal{S}(x)$ and using the fact that $ \alpha \geq 1 $ we arrive at
\[
\| T \|_v =  \|_{0} ( P_0 ) \|_{\upsilon} \leq  b_1 ( \alpha ) := \mathcal{S}^{\ast} ( \lambda (1 - \alpha ) )   - \mathcal{S}^{*} (\lambda)  ,
\]
provided that $ \alpha $ is such that 
\begin{equation}\label{eq:condition2}
\mathcal{S}^{\ast} ( \lambda (1 - \alpha ) ) < \infty  .
 \end{equation}
Furthermore, using \eqref{eq:taboo2} one obtains
\[
|| \pi_0^\top ||_v \leq b_2 ( \alpha ) := \frac{\sum_i  \pi_{0} (i) P_0 ( i , 0 )  }{ 1 - b_1 (\alpha) } = \frac{ \pi_{0} (0) }{ 1 - b_1 (\alpha) }.
\]
We now turn to computing a bound for  $  || P_1 - P_0||_v $.
For $i = 0$:
    \begin{eqnarray*}
       & &    \sum\limits_{j \geq 0} \alpha^{j} | P_1 ( 0 , j ) - P_0 ( 0 , j) | \nonumber
       \\
& & \quad        =  \sum\limits_{j \geq 0} \alpha^{j} \left|
\frac{r}{r + \lambda} \left ( \frac{\lambda}{\lambda +r } \right )^{j-1}1_{ j \geq 1 }
 - \int_{0}^{\infty} e^{-\lambda x} \frac{(\lambda x)^{j}}{ j! } d\mathcal{S}(x)\right|
\nonumber  \\
& & \quad        =
\int_0^\infty e^{-\lambda x}  d\mathcal{S}(x) +
\sum\limits_{j \geq 0} \alpha^{j+1} \left|
\frac{r}{r + \lambda} \left ( \frac{\lambda}{\lambda +r } \right )^{j}
 - \int_{0}^{\infty} e^{-\lambda x} \frac{(\lambda x)^{j+1}}{ (j+1)! } d\mathcal{S}(x)\right|
. \nonumber
    \end{eqnarray*}
For $i \geq 1$:
    \begin{eqnarray*}
 & &   \frac{ 1 }{ \alpha^{i} } \sum\limits_{j \geq 0} \alpha^{j} | P_1 ( i ,j  ) - P_0 ( i ,j) |
\nonumber \\
& & \qquad =
       \frac{ 1 }{ \alpha^{i} }
\sum\limits_{j \geq 0} \alpha^{j } \left|
\frac{r}{r + \lambda} \left ( \frac{\lambda}{\lambda +r } \right )^{j-i}1_{ j \geq i }
 - \int_{0}^{\infty} e^{-\lambda x} \frac{(\lambda x)^{j-i+1}}{ (j-i+1)! } d\mathcal{S}(x)\right| 1_{ j -i+1 \geq 0}
\nonumber\\
& & \qquad = \frac{1}{\alpha}
\int_0^\infty e^{ - \lambda x } d \mathcal{S} ( d x )  +
       \frac{ 1 }{ \alpha^{i} }
\sum\limits_{j \geq i} \alpha^{j +1} \left|
\frac{r}{r + \lambda} \left ( \frac{\lambda}{\lambda +r } \right )^{j-i}
 - \int_{0}^{\infty} e^{-\lambda x} \frac{(\lambda x)^{j-i+1}}{ (j-i+1)! } d\mathcal{S}(x)\right|
\nonumber\\
& & \qquad \leq 
\int_0^\infty e^{ - \lambda x } d \mathcal{S} ( d x )  +
  \sum\limits_{j \geq 0} \alpha^{j+1 } \left|
\frac{r}{r + \lambda} \left ( \frac{\lambda}{\lambda +r } \right )^{j}
 - \int_{0}^{\infty} e^{-\lambda x} \frac{(\lambda x)^{j+1}}{ (j+1)! } d\mathcal{S}(x)\right|  . 
 \nonumber
    \end{eqnarray*}
Combining the above results we let
\begin{eqnarray*}
 b_3 ( \alpha )  & := &  
 \int_0^\infty e^{-\lambda x}  d\mathcal{S}(x) +
\sum\limits_{j \geq 0} \alpha^{j+1} \left|
\frac{r}{r + \lambda} \left ( \frac{\lambda}{\lambda +r } \right )^{j}
 - \int_{0}^{\infty} e^{-\lambda x} \frac{(\lambda x)^{j+1}}{ (j+1)! } d\mathcal{S}(x)\right|
  \end{eqnarray*}
and obtain
\[
|| P_1 - P_0 ||_v \leq b_3 ( \alpha ) . 
\]
Inserting the above bounds into (\ref{eq:SSB}) we obtain as SSB
\begin{eqnarray*}
|| \pi_\theta^\top -  \pi_0^\top  ||_v  & \leq  &
 b_2 ( \alpha )
 \frac{ \theta  (1 +b_2 ( \alpha ) ) b_3 (\alpha)}{ 1 - b_1 ( \alpha ) -  \theta   (1 +b_2 ( \alpha ) ) b_3 (\alpha)}  ,
\end{eqnarray*}
provided  that
\[
\theta <  \frac{1 - b_1 ( \alpha ) }{  (1 +b_2 ( \alpha ) )b_3 (\alpha)} 
\]
and $ 1 \leq \alpha \leq \min ( 1/ \lambda , z_\lambda  )$, where $z_\lambda $ denotes the right point of the domain of the values for $ \alpha $ such that 
$ \mathcal{S}^{\ast} ( \lambda (1 - \alpha ) )  $ is finite (the case $ z_\lambda = \infty$ is not excluded). 

\begin{example}\label{ex:erer}
If  the service times are exponentially distributed with rate $ \mu $
it holds that 
\[
\mathcal{S}^{\ast} ( \lambda (1 - \alpha ) )  = \frac{\mu}{\mu + \lambda ( 1 - \alpha ) } 
\]
and $ z_\lambda = \frac{\mu + \lambda}{\lambda} - \epsilon $, for $ \epsilon > 0 $. 
The above bounds can now be explicitly computed: 
\[
b_1 (\alpha) =
\frac{\mu}{\mu +  \lambda (1 - \alpha ) }-  \frac{\mu}{\mu + \lambda } = \frac{\lambda \mu \alpha}{(\mu + \lambda )(\mu + \lambda (1- \alpha ))} ,
\]
\[
 b_2 ( \alpha ) = \frac{1-\lambda / \mu}{1 - b_1 ( \alpha )}.
 \]
 and 
 \[
 b_3 ( \alpha ) = \frac{\mu}{\lambda + \mu}  + \alpha
\sum\limits_{j \geq 0} \alpha^{j} \left|
\frac{r}{r + \lambda} \left ( \frac{\lambda}{\lambda +r } \right )^{j}
 - \left ( \frac{\lambda}{ \lambda + \mu } \right )^{j+1} \frac{\mu}{\mu + \lambda}\right| .
 \]

 Note that in case $ \mu = r $, $ b_3 ( \alpha ) $ simplifies to
 \[
  b_3 ( \alpha ) = \frac{\mu}{\lambda + \mu}  +
\alpha  \sum\limits_{j \geq 0} \alpha^{j} \left(\frac{\mu}{\mu + \lambda}\right)^2 \left ( \frac{\lambda}{ \lambda + \mu } \right )^{j}
= \frac{\mu}{\lambda + \mu}\left( 1 + \frac{\alpha \mu}{\mu+ \lambda - \alpha\lambda} \right)
 \]
 provided that $ \alpha < \frac{\lambda+\mu}{\lambda }$.
\end{example}

In the following, we let  $ \lambda=  0.5$, $ \mu = 1 $, $ r = 1 $ and $ f (s ) = 0 $ for $ s \leq 2 $ and $ f ( s ) = 1 $ for $ s >  2 $, i.e.,
we are interested in the probability of having more than 2 customers at the queue in stationary regime, i.e.,
\[
|| f ||_v =  \frac{1}{ \alpha^{3}} .
\]
For ease of computation we assume that the service times are exponentially distributed.

\begin{figure}[ht]
\begin{picture}(400,250)
\includegraphics[scale=0.8]{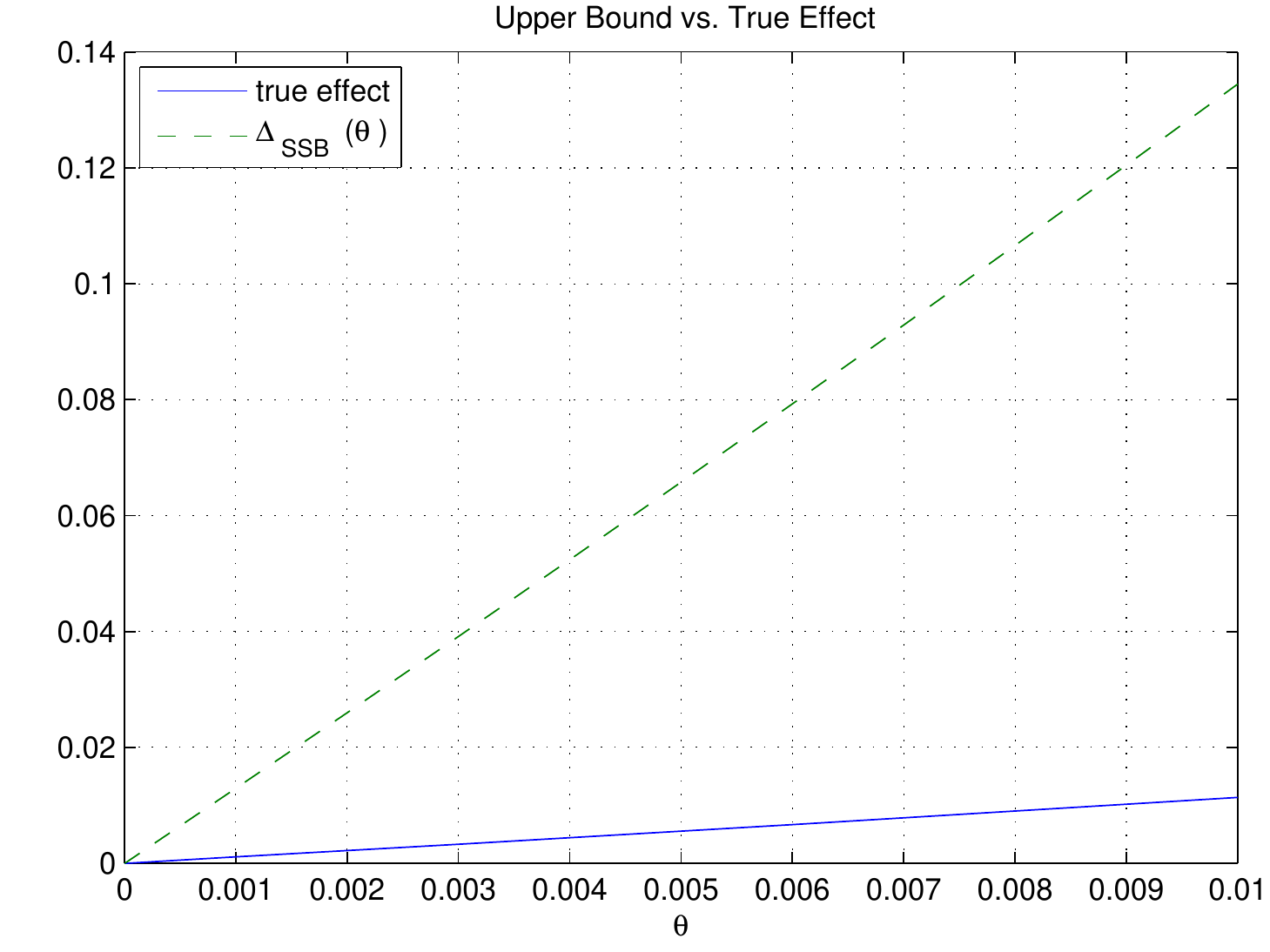}
\end{picture}
\caption{The true change in probability of more than 2 customers in the system vs. the strong stability bound.}
\label{fig:ss2}
\end{figure}

We are now able to apply the bound provided in Lemma~\ref{le:norms} to $ | \pi_\theta f - \pi_0 f | $ in combination with the above SSB, where we let $ \theta $ vary from 0 to 0.01, see Figure~\ref{fig:ss2}.
The minimization with respect to $ \alpha  $ in (\ref{eq:starrr}) has been solved numerically.
As can be seen from Figure~\ref{fig:ss2}, SSB provides qualitative insight rather than numerically satisfying approximations.

Recall that $ T = _0 \!  ( P_0 ) $ and, by Remark~\ref{re:ssb}, applicability of SSB implies stability of the system with breakdowns.  
SSB can thus be used as  means of establishing a lower bound for the domain of stability of the queue with breakdowns.
More precisely,   by  Example~\ref{ex:erer},  for $ \mu = r=1 $ condition
\[
|| T ||_v \leq b_1 ( \alpha ) < 1
\]
implies
\[
\alpha \leq \frac{ ( \mu + \lambda )^2 }{ ( 2 \mu + \lambda ) \lambda } ,
\]
which yields for the numerical setting of our example
\[
\alpha \leq \frac{9}{5} .
\]
In accordance with Corollary~\ref{cor:suff}, a lower bound for the region of stability of $ P ( \theta ) $ is
\[
\frac{1 - || T ||_v}{|| P_1 - P_0 ||_v } \geq \max_{ 1 \leq \alpha \leq 9/5 } \frac{(\mu+\lambda)^2 - \lambda(2\mu+\lambda)\alpha}{\mu(\mu+\lambda+\alpha(\mu-\lambda))},
\]
where we used the bounds provided in   Example~\ref{ex:erer}.
For the numerical values of the example we obtain
\[
\max_{ 1 \leq \alpha \leq 9/5 }  \frac{9 - 5 \alpha}{6 +   2 \alpha}   = \frac{1}{2}  ,
\]
where the maximum is attained at $ \alpha = 1 $.
Hence, the system remains stable for a breakdown probability up to $  \approx 1/2 $.

In the following section, we will show that the series expansion bound yields numerically better bounds. This comes, however, at the price of restricting the analysis to a finite version of model.

\subsection{The M/G/1/N Queue with Breakdowns  (Finite State Space)}\label{sec:A6}
In this section a M/G/1/N queue is considered with finite size $N$ (where $N$ is not too large). 
In this case the state space is $S=\{0,1,\dots,N\}$, and 
 $ D_\theta $ (short for $D_{P_\theta}$) as well as $ \pi_\theta $ (short for $\pi_{P_\theta}$) can be easily computed
numerically. 
In this case, SEB can be used for numerical computations.
We illustrate the series expansion bound  with some numerical examples.
We choose $ N =50 $ as the maximum number of jobs in the system.
Like in the previous section, we let $ \lambda =0.5$, $ \mu = 1 $,  $ r = 1 $, and assume that service times are exponentially distributed.

\begin{remark}
Note that for  large $ N $ the mean queue length of
the finite system is (almost) identical to that of the infinite one.
In this case one could use the strong stability bounds for approximate performance evaluation rather than computing SEB explicitly.
\end{remark}

We compute SEB for the $v$-norm with $\alpha = 1$. We have to check the condition put forward in (iv) of Theorem~\ref{th:relerror} numerically.
For our numerical setting we obtain
$  || ( P_{1} - P_0 ) D_0 ||_{ v }   = 8 $, which implies $  \theta  || ( P_{1} - P_0 ) D_0 ||_{ v }  < 1 $  for $ 0 \leq \theta \leq \theta_0 < 1 /8 $.
 In the following we choose $ \theta_0 = 0.1 $.

\begin{figure}[ht]
\begin{picture}(400,250)
 \includegraphics[scale=0.8]{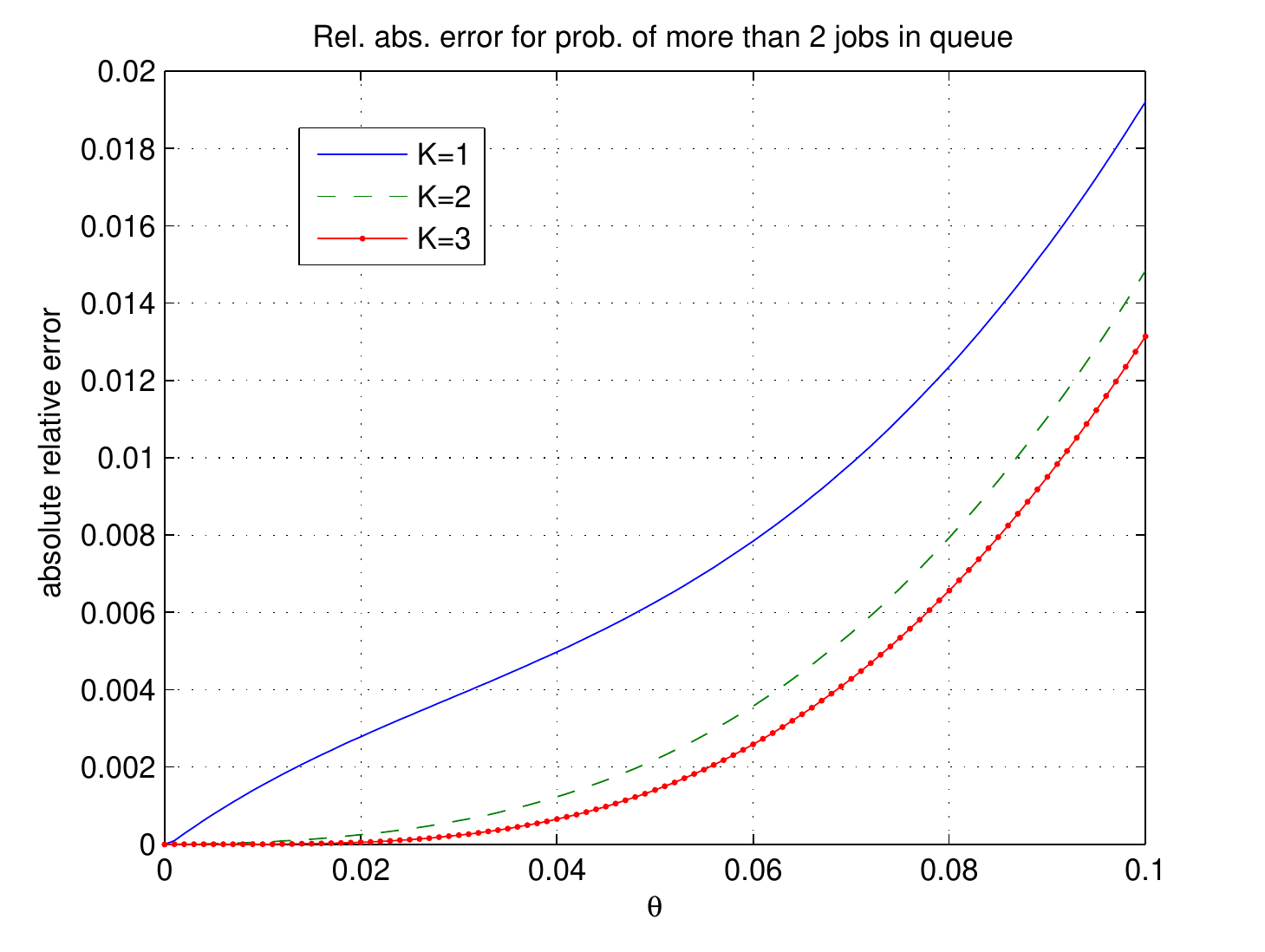}
\end{picture}
\caption{The relative absolute error for approximating the $  |  \pi_\theta^\top  f  - \pi_0^\top f | $   with SEB($K$) with $ K=1,2$ and $3$.}
\label{fig:1}
\end{figure}

In Figure~\ref{fig:1} we plot the relative absolute error of SEB($K$) for $ K=1,2 $ and $3$,
for the probability of having more than 2 customers in the systems.
More specifically, we bound $ | \pi_\theta^\top f  - \pi_0^\top f | $ for $ \theta \in [ 0 , \theta_0 ] $, with $ \theta_0 = 0.1 $,
using SEB($K$), where
$ f ( s ) = 1 $  if $ s > 2 $ and zero otherwise.
It thus holds that $ || f ||_v = 1 $.
 In line with Lemma~\ref{le:norms}, we obtain the bound  
\[
| \pi_\theta^\top f  - \pi_0^\top f | \leq \Delta_{ \rm{SEB} ( K )} ( P ( \theta ) , P_0 )  .
\]
We plot in Figure~\ref{fig:1} the absolute relative error, given by
\[
\frac{ \left| \Delta_{ \rm{SEB} ( K )} ( P ( \theta ) , P_0 ) - |  \pi_\theta^\top  f  - \pi_0^\top f | \, \right| }{ |  \pi_\theta^\top  f  - \pi_0^\top f | }  ,
\]
for $ K=1,2,3$ and $ \theta \in [ 0 , 0.1 ] $.

\subsection{Discussion of Results}
In this section  we discussed numerical approximations for the single server queue with breakdowns. SSB has the advantage of providing bounds for infinite queues, unfortunately, the numerical quality of the bounds is rather poor.
In light of Theorem~\ref{th:relerror}, this comes as no surprise.
SEB proved to be numerically very efficient for the model but required that a finite queue is studied.
There is, however, an interesting link between the two approaches as the techniques developed for SSB lend themselves to establish lower bounds of convergence for series expansions.

\section{Conclusion}
Perturbation bounds for Markov chains have been intensively studied in the literature.
Condition number bounds are attractive as they provide uniform perturbation bounds. 
Unfortunately, due to their simple structure they fail to capture
the true non-linear dependence of the stationary distribution on the Markov kernel. 
SSB, which provides a non-linear expression in the size of the perturbation, overcomes this drawback and is the only bound applicable in case 
of an infinite state space.  
We introduced a new family of bounds based on a series expansion approach.
As illustrated by a series of examples both analytical and numerical, our new bounds 
yield good results and have the desirable property that the relative error vanishes when the size of the perturbation tends to zero.
A realistic example from queueing theory illustrated the potential use of perturbation bounds in robustness analysis.

\section*{Acknowledgement}
The authors are grateful to an anonymous reviewer for valuable remarks on an earlier version of the paper.



\end{document}